\documentstyle{article}

\author{
        {\bf Ladislav HLAVAT\'Y} and {\bf Libor \v{S}NOBL}\thanks{e-mail  :hlavaty@br.fjfi.cvut.cz, snobl@gauss.fjfi.cvut.cz} \\  
        {\it  Department of Physics, Faculty of Nuclear Science,}\\
        {\it Czech Technical University }\\
        {\it Brehova 7, 115 19 Prague 1, Czech Republic}\\
        }

\title{Solution of a Yang--Baxter system}

\def \be {\begin{equation}}
\def \bea {\begin{eqnarray}}
\def \ee {\end{equation}}
\def \eea {\end{eqnarray}}

\def \ll {\label}

\def \ox {\otimes}
\def \unit {{\bf 1}}
\def \complex {{\bf C}}

\def \newblock {}
\def \rf {(\ref}
\def \eqn {equation}
\def \sln {solution}
\def \YB {Yang--Baxter}
\def \YBE {Yang--Baxter equation}
\def \YBS {Yang--Baxter system}
\def \cfn {classification}
\def \reln {relation}
\def \tfn {transformation}
\def \mt {matri}
\begin{document}
\bibliographystyle{unsrt}    
                             
\maketitle

\begin{abstract}
Yang--Baxter system related to quantum doubles is introduced and large class of both continuous and discrete symmetries of the solution manifold are determined.
Strategy for solution of the system  based on the symmetries  is suggested and accomplished in the dimension two.
 The complete list of invertible \sln s of the system is presented.
\end{abstract}

\section{Introduction}
The Yang--Baxter equations, both constant and spectral dependent proved to be an important tool for various branches of theoretical physics. They represent a system of $N^6$ cubic \eqn s for elements of $N^2\times N^2$ matrix $R$ and can be written in the well known form
\begin{equation} R_{12}R_{13}R_{23}     =
R_{23}R_{13}R_{12}.
\ll{cybe}\end{equation}
Even though many \sln s are known for all types of the \YBE s \cite{sogo:ybecfn,hla:usybe,hla:baxtn,sww:cfn6v,wang:cfn8v},  
until now the complete solution is known only for the constant \YBE s in the dimension two \cite{hie:ybecfn}, i.e. matrices $4\times 4$.

Various extensions of the \YBE s for several matrices, called \YBS s,
appeared in literature \cite{frimai,nijetal:pra92,hla:qbg}.  The constant systems are used mainly for construction of special Hopf algebras while the  spectral dependent \sln s are applied in quantum integrable models.
Examples of both types together with their 
particular solutions were presented in \cite{hla:ybssa}. 

As the \YBS s usually contain several \YB--type \eqn s, it is convenient to introduce the following notation: 
{\em The Yang--Baxter commutator} $[R,S,T]$ of (constant) $N^2\times N^2$ matrices $R,S,T$ is $N^3\times N^3$ matrix
\be [R,S,T]:=R_{12}S_{13}T_{23}-
T_{23}S_{13}R_{12}. \ll{cybc} \ee
In this notation, the constant \YBE {} is written as
\be [R,R,R]=0,\ \ee
while e.g. the system for quantized braided groups 
\cite{hla:qbg} reads
\begin{equation} [Q,Q,Q]=0,\ [R,R,R]=0,\label{rrr0}
\end{equation}
\begin{equation}
[Q,R,R]=0,\ [R,R,Q]=0. \label{zzr0}
\end{equation}
The complete set of invertible solutions of this system in the dimension two is given in \cite{hla:cfnqbg}.

The goal of the present paper is to give a complete \sln {} of another, more complicated constant system related to quantum doubles i.e. special quasitriangular Hopf algebras constructed from the tensor product of Hopf algebras by defining a pairing between them.

\section{WXZ system and its symmetries}
In the paper  \cite{vlad:qdouble} a method of obtaining the quantum doubles for  pairs of FRT quantum groups is presented: 
Let two quantum groups are given by relations \cite{FRT}
\[ W_{12}U_{1}U_{2}=U_{2}U_{1}W_{12} \]
\[ Z_{12}T_{1}T_{2}=T_{2}T_{1}Z_{12} \]
where $W$ and $Z$ are matrices $N^2\times N^2$ satisfying the \YBE s
\begin{equation} [W,W,W]=0,\ll{www}\ee
\be [Z,Z,Z]=0, \ll{zzz}\ee
and suppose that there is a \mt x $X$ that satisfies the \eqn s
\begin{equation} [W,X,X]=0,\ll{wxx1}\ee
\be [X,X,Z]=0, \ll{xxz}\end{equation}
Then the \reln s
\be X_{12}U_{1}T_{2}=T_{2}U_{1}X_{12} \ll{wxzs}\ee
define  quantum double with the pairing
\be <U_{1},T_{2}>=X_{12}. \ee
On the other hand, the \eqn s \rf{www})--\rf{xxz}), that in the following we shall call the {\em WXZ system}, are constant version of the spectral dependent \YBS s for nonultralocal models presented in \cite{hlakun:nulm}.

There are two natural questions related to the WXZ system:
\begin {itemize}
\item
Is there a \mt x $X$ such that for any pair of matrices $W,Z$ that solve the \YBE s  the triple $(W,X,Z)$ solves the WXZ system? 
\item
Is there for any \mt x $X$ a pair of matrices $W,Z$  such that the triple $(W,X,Z)$ solves the WXZ system?
\end {itemize}
Answers to both these questions are positive because  the following two simple propositions hold
\begin {itemize}
\item
Let  $W,Z$ are arbitrary\ solutions of the \YBE s
$[W,W,W]=0,\ [Z,Z,Z]=0$. Then the triple  
$(W, X=\unit,Z)$ is a solution of the  system \rf{www})--\rf{xxz}).
\item
Let $X$ is an arbitrary  \mt x $N^2\times N^2$ and  $P$ is the permutation matrix. Then the triple  
$(W=P, X,Z=P)$ is a solution of the  system \rf{www})--\rf{xxz}).
\end {itemize}
Other simple solution of the WXZ system is $(W=R,X=R,Z=R)$, where $R$ is an arbitrary \sln{} of the \YBE.

Besides the above mentioned, one can find  solutions of the WXZ system from the knowledge of solutions of the system (\ref{rrr0})--(\ref{zzr0}). Namely, if $(Q,\, R)$ is a solution of
the system of Yang--Baxter type equations (\ref{rrr0})--(\ref{zzr0})
then $(W=Q,X=R,Z=Q)$ and $(W=Q,X=R,Z=R^+QR^-)$ are solutions of the
system \rf{www})--\rf{xxz}). 

Solution of the system is essentially facilitated by knowledge of its symmetries. It is easy to check that
the set of \sln s is invariant under both continuous transformations
\[ W'=\omega(T\ox T)W(T\ox T)^{-1}\]
\be X'=\xi(T\ox S)X(T\ox S)^{-1}\ll{wxzsym}\ee
\[ Z'=\zeta(S\ox S)Z(S\ox S)^{-1}\]
where
\[ \omega,\xi,\zeta\in\complex,\ \ T,S\in {\cal GL}(N,\complex), \]
and discrete transformations
\be (W',X',Z')=(W^T,X^T,Z^T) \ll{dsym}\ee
\be (W',X',Z')=(W^a,X,Z^b), \ \ a=id,\# ,\ b=id,\#  \ll{dsym1}\ee
\be (W',X',Z')=(W^c,X^-,Z^d), \ \ c=+,-,\ d=+,- \ll{dsym2}\ee
\be (W',X',Z')=(Z^c,X^+,W^d), \ \ c=+,-,\ d=+,- \ll{dsym3}\ee
where $Y^T$ is transpose of $Y$, $Y^+:=PYP,\ P$ being the permutation \mt x, $Y^-:=Y^{-1},$ $Y^{id}:=Y,\ Y^\# :=(Y^+)^-=(Y^-)^+$.

These symmetries can be of course composed, so that one can define e.g. the "antidiagonal transposition"
\be Y^{at}:=(\sigma\otimes\sigma) Y^T (\sigma\otimes\sigma) \ 
{\rm where} \ \sigma=\left( \begin{array}{cc} 0&1\\1&0 \end{array} \right)   
\ll{at}\ee
and the symmetry
\be(W',X',Z')=(W^{at},X^{at},Z^{at}) \ll{atsym}\ee

Beside that there are useful  "conditional symmetries"  that can be expressed as

{\bf Lemma 1}: Let $W,X,Z$ are  \mt ces $N^2\times N^2$ that  solve WXZ system \rf{www})--\rf{xxz}) and  $T,A$ are \mt ces  $N\times N$ such that 
\be [W,T\ox T]=0,\ee
\be X(T\ox\unit)=(A\ox\unit)X\ {\rm or}\ (T\ox\unit)X=X(A\ox\unit).\ee
Then 
\be  W'=W,\  X'=(T\ox \unit) X ,\  Z'=Z\ll{condsym01}\ee 
\be  W'=W,\  X'=X(T\ox \unit)  ,\  Z'=Z\ll{condsym11}\ee 
and
\be  W'=(T\ox \unit)W(T\ox \unit)^{-1},\  X'=X,\  Z'=Z\ll{condsym1}\ee 
solve the WXZ system  as well.

{\bf Lemma 2}: Let $W,X,Z$ are  \mt ces $N^2\times N^2$ that  solve WXZ system \rf{www})--\rf{xxz}) and  $S,A$ are \mt ces  $N\times N$ such that 
\be [Z,S\ox S]=0,\ee
\be X(S\ox\unit)=(A\ox\unit)X\ {\rm or}\ (S\ox\unit)X=X(A\ox\unit).\ee
Then 
\be  W'=W,\  X'=(S\ox \unit) X ,\  Z'=Z\ll{condsym02}\ee 
\be  W'=W,\  X'=X(S\ox \unit)  ,\  Z'=Z\ll{condsym12}\ee 
and
\be  W'=W,\  X'=X,\  Z'=(S\ox \unit)Z(S\ox \unit)^{-1}\ll{condsym2}\ee 
solve the WXZ system  as well.

Proofs of these lemmas can be done by direct calculations.
There are other symmetries of the WXZ system that are extensions of the twisting \tfn s of the solutions  \YBE {} but we are not going to use them in the following as it seems that they do not produce equivalent quantum doubles.

\section{Solution of the WXZ system in the dimension two}\ll{swxzind2}
Even for the lowest nontrivial dimension two, solution of the system \rf{www})--\rf{xxz}) represents a tremendous task, namely solving 256 cubic equations for 48 unknowns. That's why it is understandable that the assistance of computer programs for symbolic calculations is essential in the following\footnote{Reduce 3.6 and Maple V was used}. On the other hand it does not mean that one can find the solutions by pure brute force, namely applying a procedure SOLVE to the system of the 256 equations. Moreover, the interpretation of the result would be extremely difficult as  hundreds of solutions would appear, many of them related by the symmetries.

There are several strategies for solution of the problem given above. All of them are based on the symmetries of the system and the knowledge of the complete set of solutions of the \YBE{} in the dimension two \cite{hie:ybecfn}.

One possible (and obvious) strategy is solving the \eqn s \rf{wxx1}),\rf{xxz}) for all pairs of matrices $W,Z$ in the Hietarinta's list of \sln s of the \YBE{}. By this way we reduce the problem to 128 quadratic \eqn s for 16-22 unknowns (depending on the number of parameters in the \sln s of the \YBE{}).

Another strategy is to use the symmetry \rf{wxzsym}) to simplify the matrix $X$ as much as possible, then solve the (linear in $W$ or $Z$) \eqn s \rf{wxx1}),\rf{xxz}) and finally solve the \YBE s \rf{www}),\rf{zzz})  for $W$ and $Z$. This is an analogue of the strategy accepted in \cite{hla:cfnqbg} for solving the system \rf{rrr0})--\rf{zzr0}).

Both the mentioned strategies yield sets of equations and unknowns that are still too large to be solved by  the computer program.
The strategy that seems to work is a combination of the two previous. 

First we find the matrices $X$ that solve the equation \rf{wxx1}) for each invertible $W$ that belongs to the Hietarinta's list.
In the next step we solve the \eqn{}  \rf{xxz}) 

and finally determine $Z$ from \rf{zzz}) using the results of the previous step. 
As the list of \sln s is quite large it is essential in each step to factorize its results by the symmetries \rf{wxzsym})--\rf{dsym3}).

Detailed description of this procedure and results are presented in the following subsections.

\subsection{Solutions of the \YBE s}
As it was mentioned, the important point is the knowledge of all invertible \sln s of the constant \YBE{} $[W,W,W]=0$. Their list up to the symmetries, denoted  ${\cal S}$, 
consists of twelve (intersecting) subsets \cite{hie:ybecfn,hla:usybe,hla:baxtn} parametrized by up to three complex numbers. 
\[ {\cal S}:=\{R_{3.1}(r,s,t),R_{2.1}(r,s),R_{2.2}(r,s),R_{2.3}(u,v), R_{1.1}(r),R_{1.2}(s),R_{1.3}(u), \]
\be R_{1.4}(t),R_{0.1}, R_{0.2},R_{0.3},
R_{0.4} \ 
|\ r,s,t\in\complex\setminus\{0\} ,\ u,v\in\complex\} 
\ll{ybesms}\ee
The matrices $R$ are defined in the Appendix \ref{apybeslns}.

For later use we shall define three subsets of ${\cal S}$, namely
\be {\cal S}_{5V}:=\{R_{2.1}(r,s),R_{2.2}(r,s), R_{3.1}(r,s,t)|\ r,s,t\in\complex\setminus\{0\}\} ,
\ll{6vms}\ee
\[ {\cal S}_{8V}:=\{R_{3.1}(r,s,t),R_{2.1}(r,s),R_{2.2}(r,s),R_{1.1}(r),R_{1.2}(s),R_{1.4}(t), \]
\be R_{0.1}, R_{0.2},R_{0.3},
R_{0.4}| \ r,s,t\in\complex\setminus\{0\} \ \} 
\ll{8vms}\ee
\be {\cal S}_{ST}:=\{R_{3.1}(1,1,1),R_{2.3}(u,v),R_{1.3}(u),R_{0.1}|\ u,v\in\complex\} .
\ll{stms}\ee
They are the sets of \sln s of the "five--vertex" form 
\be R_{5V}  =  \left (\begin {array}{cccc} a&0&0&0\\\noalign{\medskip}0&b&0
&0\\\noalign{\medskip}0&q&c&0\\\noalign{\medskip}0&0&0&d\end {array}
\right ), \ee
"eight--vertex" form
\be R_{8V}  = \left (\begin {array}{cccc} a&0&0&x\\\noalign{\medskip}0&b&p
&0\\\noalign{\medskip}0&q&c&0\\\noalign{\medskip}y&0&0&d\end {array}
\right ), \ee
and "special triangular" form 
\be R_{ST}  =  \left (\begin {array}{cccc} 1&0&0&0\\\noalign{\medskip}a&1&0
&0\\\noalign{\medskip}b&0&1&0\\\noalign{\medskip}c&d&e&1\end {array}
\right ). \ee
It is easy to check that 
\be {\cal S}={\cal S}_{8V}\cup {\cal S}_{ST}, \ {\cal S}_{5V}\subset{\cal S}_{8V} \ee
The reason for this partition is that these subsets have common symmetries and common \mt ces $X$ that solve the \eqn{} \rf{wxx1}).
Namely, each element of  ${\cal S}_{5V}$ is invariant under 
\be R'=(A\ox A) R (A\ox A)^{-1} \ll{rtr}\ee where
\be A=\left (\begin {array}{cc} a&0\\0&b \end {array}\right ) \ll{t6v}\ee
while each element of  ${\cal S}_{8V}$ is invariant under \rf{rtr}) where
\be A=\left (\begin {array}{cc} a&0\\0&-a \end {array}\right ) \ll{t8v}\ee
and each element of  ${\cal S}_{ST}$ is invariant under \rf{rtr})  where
\be A=\left (\begin {array}{cc} a&0\\b&a \end {array}\right ). \ll{tst}\ee
These invariances imply the existence of solutions $X=A\ox B$ of the  equation \rf{wxx1}) common to all elements of the subsets ${\cal S}_{5V},\ {\cal S}_{8V}$ or ${\cal S}_{ST}$. Actually we shall find that slightly more general solutions are generic for each of the subsets.

\subsection{Solutions of the WXX equation}\ll{wxxsolns}
The goal of this subsection is to find all invertible \mt ces $X$ that solve the \eqn{} 
\be [W,X,X] =0\ll{wxx}\ee
for all invertible \sln s $W$ of the \YBE {} \rf{www}). 
First of all we can exploit the symmetry \rf{wxzsym}) to reduce the \mt x $X$ to one of the forms $A_{1},\ldots,A_{14}$ given in the Appendix \ref{xansatz}.
Then we choose a solution $W$ of the \YBE{} and solve $X$ from the \eqn{} \rf{wxx}).
The list of the \sln s up to the WXZ symmetries \rf{wxzsym}), \rf{dsym1}),\rf{dsym3}), \rf{condsym01}), \rf{condsym11})
can be described in terms of 23 classes of \mt ces 
parametrized by complex numbers and
given in the Appendix \ref{xseznam}.

As mentioned in the previous subsection, there are \sln s common to all elements in subsets ${\cal S}_{5V},\ {\cal S}_{8V}$ or ${\cal S}_{ST}$. Namely:
\begin{itemize}
\item If $W\in{\cal S}_{5V}$ then $X_{\ref{diag}}(a,b,c)$ and $X_{\ref{X.jj1}}(a,b,c)$ solve \rf{wxx}) for all $a,b,c\in\complex$. 
\item If $W\in{\cal S}_{8V}$ then $X_{\ref{diag}}(a,\epsilon_{1},\epsilon_{2} a ),$   and $X_{\ref{X.jj1}}(\epsilon_{1},a,\epsilon_{1}a)$ solve \rf{wxx})
for all $ a\in \complex$ and $\epsilon_{1},\epsilon_{2}=\pm 1$.
\item If $W\in{\cal S}_{ST}$ then $X_{\ref{X.jj2}}(a,b,c)$ and $X_{\ref{X10}}(a,b,c)$ solve \rf{wxx}) for all $ a,b,c\in\complex$.
\end {itemize}
We shall call these \sln s {\em generic}. Notation of the X--\mt ces above as well as below corresponds to that  in the Appendix \ref{xseznam}.

Below we display all invertible non--generic \sln s of \rf{wxx}).
The list 
is complete up to the WXZ symmetries \rf{wxzsym}), \rf{dsym1}),\rf{dsym3}), \rf{condsym01}), \rf{condsym11}). That's why we also present the \tfn s that leave given $W$ invariant or form invariant, i.e. invariant up to the change of parameters.

The ranges of parameters are
\[r,s,t\in\complex\setminus\{ 0 \}, \ \ 
\ u,v,a,b,c,d,e,f,g,h\in\complex, \ \ 
\ \epsilon,\epsilon_{1},\epsilon_{2}=\pm 1. \]
\\{\bf List of non--generic solutions of the equation \rf{wxx}):}
\begin{itemize} 

\item $W=R_{3.1}(r,s,t)$.
\\
This \mt x is invariant w.r.t. \rf{rtr}),\rf{t6v}) and
\be  R_{3.1}(r,s,t)^{T}=R_{3.1}(r,s,t),\ R_{3.1}(r,s,t)^{-1}=R_{3.1}(r^{-1},s^{-1},t^{-1}),\ee 
\be R_{3.1}(r,s,t)^{+}=R_{3.1}(s,r,t) \ee
\be(T\ox T)R_{3.1}(r,s,t)(T\ox T)^{-1}=t\,R_{3.1}(\frac{s}{t}, \frac{r}{t},\frac{1}{t}),\ {\rm where }\ T=\left (\begin {array}{cc} 0&1\\1&0 \end {array}\right ). 
\ee
The  non--generic \sln s of the \eqn {} $[W,X,X]=0$ up to the symmetries exist only for special values of the parameters $r,s,t$.
     \begin{enumerate} 
     \item $\ \ W=R_{3.1}(s,s,1),\ X=X_{\ref{X8.3}}(a,b,c) 
     $
     \item $\ \ W=R_{3.1}(s,s,1),\ X=X_{\ref{X8.4}}(a,b,c) 
     $
     \item $\ \ W=R_{3.1}(s,s,1),\ X=X_{\ref{X8.5}}(a,b) 
     $
     \item $\ \ W=R_{3.1}(s,-s,1),\quad X=X_{\ref{X8.7}}(a) 
     $
     \item $\ \ W=R_{3.1}(1/s,s,\epsilon), \quad X=X_{\ref{h5}}(s,a,\epsilon sa) 
            $
     \item $\ \ W=R_{3.1}(\epsilon,\epsilon,-\epsilon),  \quad X=X_{\ref{h6}}(\epsilon,b,-b,c) $
     \item $\ \ W=R_{3.1}(-1,-1,1), \quad X=X_{\ref{X.30}}(a,b,c,d) 
     $
     \item $\ \ W=R_{3.1}(1,1,1), \quad X=X_{\ref{X8.15}}(a,b,c,d,e,f,g,h) $
     \item $\ \ W=R_{3.1}(1,1,1), \quad X=X_{\ref{X8.16}}(a,b,c,d,e,f,g,h) $
     \end{enumerate} 
 
\item $W=R_{2.1}(r,s)$ \\
This \mt x is invariant w.r.t. \rf{rtr}),\rf{t6v}) and
\be  R_{2.1}^{at}(r,s)=R_{2.1}(s,r),\ R_{2.1}(r,s)^{-1}=R_{2.1}(r^{-1},s^{-1}) \ee
\be(T\ox T)R_{2.1}^{+}(r,s)(T\ox T)^{-1}=R_{2.1}(r,s),\ {\rm where }\ T=\left (\begin {array}{cc} 0&1\\1&0 \end {array}\right ). 
\ee
The  non--generic \sln s of the \eqn {} $[W,X,X]=0$ up to the WXZ symmetries are
     \begin{enumerate} 
     \item $\ \ W=R_{2.1}(r,s),\ X=X_{\ref{h5}}(r^{-1},a,sa) 
            $
     \item $\ \ W=R_{2.1}(s,s),\ X=X_{\ref{h6}}(s^{-1},a,sa,b) 
            $
     \item $\ \ W=R_{2.1}(1,-1),\ X=X_{\ref{h7}}(a,b) 
            $
     \item $\ \ W=R_{2.1}(i,i),\ X=X_{\ref{octag}}(a) 
            $
     \item $\ \ W=R_{2.1}(i,-i), \quad X=X_{\ref{X8.7}}(a)$ (coincide with the case 4 in  $W=R_{3.1}$).
     \item $\ \ W=R_{2.1}(\epsilon,\epsilon)$, coincide with the cases 1--3,5,7--9 in  $W=R_{3.1}$.
      \end{enumerate} 
\item $W=R_{2.2}(r,s).$ 
\\
This \mt x is invariant w.r.t. \rf{rtr}),\rf{t6v}) and
\be  R_{2.2}^{at}(r,s)=-rsR_{2.2}(-r^{-1},-s^{-1}),\ R_{2.2}(r,s)^{-1}=R_{2.2}(r^{-1},s^{-1}) \ee
\be(T\ox T)R_{2.2}^{+}(r,s)(T\ox T)^{-1}=-rsR_{2.2}(-s^{-1},-r^{-1}),\ 
 {\rm where }\ T=\left (\begin {array}{cc} 0&1\\1&0 \end {array}\right ). 
\ee
The  non--generic \sln s of the \eqn {} $[W,X,X]=0$ up to the WXZ symmetries are
     \begin{enumerate} 
     \item $\ \ W=R_{2.2}(r,s),\ X=X_{\ref{h5}}(r^{-1},a,-ar^{-1}) 
            $
     \item $\ \ W=R_{2.2}(s,s),\ X=X_{\ref{h6}}(s^{-1},a,-s^{-1}a,g) 
            $
     \item $\ \ W=R_{2.2}(1,-1),\ R_{2.2}(i,i)$ coincide with cases  3,4 in $W=R_{2.1}$.
      \end{enumerate} 

\item $W=R_{2.3}(u,v)$ \\
This \mt x is invariant w.r.t. \rf{rtr}), \rf{tst}) and
\be R_{2.3}^{at}(u,v)=R_{2.3}(u,v)\ee
\be (T\ox T)R_{2.3}^{-1}(u,v)(T\ox T)^{-1}=R_{2.3}(u,2u-v),\ {\rm where }\ T=\left (\begin {array}{cc} 1&0\\0&-1 \end {array}\right ) 
\ee
\be (T\ox T)R_{2.3}^{+}(u,v)(T\ox T)^{-1}=R_{2.3}(1/u,v/u^2),\ {\rm where }\ T=\left (\begin {array}{cc} u&0\\0&1 \end {array}\right ) 
\ee
and there are no non--generic \sln s of the \eqn {} $[W,X,X]=0$ up to the WXZ symmetries.

\item $W=R_{1.1}(r).$ 
\\
This \mt x is invariant w.r.t. \rf{rtr}),\rf{t8v}) and
\be  R_{1.1}^{t}(r)=R_{1.1}(r),\ R_{1.1}(r)^{+}=R_{1.1}(r) \ee
\be(T\ox T)R_{1.1}^{-1}(r)(T\ox T)^{-1}=-\frac{1}{4}R_{1.1}(-r),\ 
 {\rm where }\ T=\left (\begin {array}{cc} 1&0\\0&i \end {array}\right ). 
\ee
The  non--generic \sln s of the \eqn {} $[W,X,X]=0$ up to the WXZ symmetries are
     \begin{enumerate} 
     \item $\ \ W=R_{1.1}(r),\ X=X_{\ref{h6}}(r,-\epsilon,\epsilon r,\epsilon) 
            $
     \item $\ \ W=R_{1.1}(r),\ X=X_{\ref{X2.1}}(a,\epsilon,\frac{r+1}{r-1}),\  r\neq 0,\pm 1$ 
     \item $\ \ W=R_{1.1}(i),\ X=X_{\ref{X9.9}}(a)\ $ 
     \end{enumerate} 

\item $W=R_{1.2}(s)$ \\
This \mt x is invariant w.r.t. \rf{rtr}), \rf{t8v}) and
\be (T\ox T)R_{1.2}^{-1}(s)(T\ox T)^{-1}=R_{1.2}(1/s),\ {\rm where }\ T=\left (\begin {array}{cc} 1&0\\0&\sqrt{s} \end {array}\right ) 
\ee
and there are the following non--generic \sln s of the \eqn {} $[W,X,X]=0$ up to the symmetries:
     \begin{enumerate}
     \item $W=R_{1.2}(s),\ X=X_{\ref{h5}}(1,\epsilon,-\epsilon)
     $
     \item  $W=R_{1.2}(s),\ X=X_{\ref{h6}}(s^{-1},\epsilon,-\epsilon s^{-1},\frac{\epsilon}{s-1}),\ $
     \item $W=R_{1.2}(s),\ X=X_{\ref{h7}}(a,b), 
     \ a^2-b^2=s+1$
     \end{enumerate}

\item $W=R_{1.3}(u)$ \\
  This \mt x is invariant w.r.t. \rf{rtr}), \rf{tst}) and
    \be  R_{1.3}^+(u)=R_{1.3}^{-1}(u)=(T\ox T)R_{1.3}(u)(T\ox T)^{-1},
  \ {\rm where }\  T=\left (\begin {array}{cc} 1&0\\0&-1 \end {array}\right ) 
  \ee
  and there are the following non--generic \sln s of the \eqn {} 
  $[W,X,X]=0$ up to the symmetries:
     \begin{enumerate}
     \item $W=R_{1.3}(u),\ X=X_{\ref{h11.2}}(a,b,b-u-1,-ua),
          \ u\neq -1$
     \item $W=R_{1.3}(1),\ X=X_{\ref{h11.1}}(a,b)$
     \end{enumerate}
\item $W=R_{1.4}(t)$.
\\
This \mt x is invariant w.r.t. \rf{rtr}),\rf{t6v}) and
  \be  R_{1.4}^{t}(t)=R_{1.4}(t),\ R_{1.4}^+=R_{1.4}(t),\ R_{1.4}^{-1}(t)=R_{1.4}(t^{-1}),\ee
\be(T\ox T)R_{1.4}(t)(T\ox T)^{-1}=R_{1.4}(t),\ {\rm where }\ T=\left (\begin {array}{cc} 0&1\\1&0 \end {array}\right ). 
\ee
The  non--generic \sln s of the \eqn {} $[W,X,X]=0$ up to the symmetries are
     \begin{enumerate} 
     \item $\ \ W=R_{1.4}(t),\ X=X_{\ref{X8.3}}(a,\epsilon,a) 
     $
     \item $\ \ W=R_{1.4}(t),\ X=X_{\ref{X8.4}}(a,1,-a) $
     \item $\ \ W=R_{1.4}(t), \quad X=X_{\ref{X9.9}}(a) 
     $
     \item $\ \ W=R_{1.4}(t),\ X=X_{\ref{X9.4}}(a) 
     $
\item  $\ \ W=R_{1.4}(\pm 1)$. This matrices can be transformed to $R_{3.1}(-1,-1,1)$ by symmetry transformations so that these solutions are equivalent to the cases 1,2,3,5,7,8 $(s=-1,\epsilon=-1)$ of $W=R_{3.1}$.
     \end{enumerate} 
\item $W=R_{0.1}$ \\
  This \mt x is invariant w.r.t. \rf{rtr}), \rf{t8v}) and
    \be  R_{0.1}^{at}=R_{0.1},\ R_{0.1}^+=R_{0.1},\ee
  \be (T\ox T)R_{0.1}^{-1}(T\ox T)^{-1}=R_{0.1},
  \ {\rm where }\  T=\left (\begin {array}{cc} 1&0\\0&i \end {array}\right ) 
  \ee
  and there are the following non--generic \sln s of the \eqn {} 
  $[W,X,X]=0$ up to the symmetries:
     \begin{enumerate}
     \item $W=R_{0.1},\ X=X_{\ref{X10.5}}(a,b,c,\epsilon)
     $
     \item $W=R_{0.1},\ X=X_{\ref{X10.6}}(a,b,c) 
     $
     \end{enumerate}
\item $W=R_{0.2}$ \\
  This \mt x is invariant w.r.t. \rf{rtr}), \rf{t8v}) and
    \be  R_{0.2}^{at}=R_{0.2},\ R_{0.2}^+=R_{0.2},\ee
  \be (T\ox T)R_{0.2}^{-1}(T\ox T)^{-1}=R_{0.2},
  \ {\rm where }\  T=\left (\begin {array}{cc} 1&0\\0&i \end {array}\right ) 
  \ee
  and there are the following non--generic \sln s of the \eqn {} 
  $[W,X,X]=0$ up to the symmetries:
     \begin{enumerate}
     \item $W=R_{0.2},\ X=X_{\ref{h5}}(-1,\epsilon,-\epsilon)$
     \item $W=R_{0.2},\ X=X_{\ref{X.30}}(\epsilon,b,0,1) 
     $
     \end{enumerate}
\item $W=R_{0.3}$ \\
  This \mt x is invariant w.r.t. \rf{rtr}), \rf{t8v}) and
    \be  R_{0.3}^T=R_{0.3},\ (T\ox T)R_{0.3}^{+}(T\ox T)^{-1}=R_{0.3},
  \ {\rm where }\  T=\left (\begin {array}{cc} 0&1\\1&0 \end {array}\right ) ,\ee
  \be 2(T\ox T)R_{0.3}^{-1}(T\ox T)^{-1}=R_{0.3},
  \ {\rm where }\  T=\left (\begin {array}{cc} 1&0\\0&i \end {array}\right ) 
  \ee
  and there are the following non--generic \sln s of the \eqn {} 
  $[W,X,X]=0$ up to the symmetries:
     \begin{enumerate}
     \item $W=R_{0.3},\ X=X_{\ref{h6}}(a,-1,a,i) $
     \item $W=R_{0.3},\ X=X_{\ref{X1.4}}(a,b) $
     \end{enumerate}
\item $W=R_{0.4}=P$ \\Arbitrary matrix X solves the equation  $[W,X,X]=0$.
\end{itemize} 

\subsection{Solution of  XXZ and ZZZ \eqn s for given $X$}\label{xxzzsolns}
\def \nogensol {\\ There are no nontrivial solutions for general values of parameters. For special values we get
   \begin{enumerate}
}

In the preceding subsection we have found all pairs $(W,X)$ of invertible $4\times 4$ matrices that solve the equations \rf{www}), \rf{wxx1}). The solutions are given in terms of the matrices $X$ presented in the Appendix \ref{xseznam}. The last step to do for solving the WXZ system is to find the invertible $4\times 4$ matrices $Z$ that for given $X$ solve the equations
\be [X,X,Z]=0,\ [Z,Z,Z]=0. \ll{xxzz}\ee

{\em Trivial solution} for any $X$ is $Z=R_{0.4}=P=$ permutation matrix. Below we display the list of nontrivial solutions of the equations \rf{xxzz}) for all $X$--matrices from the subsection \ref{xseznam}.

Let us note that in spite of the fact that the matrices $Z$ must solve the \YBE{}, they need not belong to the list of solutions in the Appendix \ref{apybeslns} because in general we have not at disposal all the \tfn s \rf{ybesyms}) up to which that list is complete, because they might already be used as a part of WXZ symmetries \rf{wxzsym}), \rf{dsym1}), \rf{dsym3}) to fix the matrices $W$ and $X$. The only symmetry that we can use in general is 
\be (W',X',Z')=(W,X,Z^\#). \ll{zsym}\ee
Beside this we can use the $WXZ$ symmetries that leave $X$ (and corresponding $W$) invariant. 
For some matrices $X$ the whole sets of the solutions of the \YBE{} ${\cal S}$ or ${\cal S}_{5V}$ or ${\cal S}_{8V}$ or ${\cal S}_{ST}$ solve the equations \rf{xxzz}). Due to the restricted symmetries of matrices $Z$ 

it is convenient for \cfn{} of solutions to introduce  another special set of \YB{} solutions 
\[ {\cal R}_{8V}:=\{R_{3.1}(r,s,t),R_{2.1}(r,s),R_{2.2}(r,s), QR_{1.1}(r)Q^{-1},\]
\[QR_{1.2}(s)Q^{-1},Q{R_{1.2}(s)}^TQ^{-1}, QR_{1.4}(t)Q^{-1},\]
\[QR_{0.1}Q^{-1},Q R_{0.2}Q^{-1},QR_{0.3}Q^{-1}, 
QR_{0.2}^TQ^{-1},QR_{0.3}^TQ^{-1},\]
\be\left ( \begin{array}{cccc} 1 & 0 & 0 & p/q^2 \\
                            0 & 1 & p & 0 \\
                            0 & p & 1 & 0 \\
                            pq^2 & 0 & 0 & 1 \end{array} \right ), \ \  
\left ( \begin{array}{cccc} 1 & 0 & 0 & p/q^2 \\
                            0 & -1 & p & 0 \\
                            0 & p & -1 & 0 \\
                            pq^2 & 0 & 0 & 1 \end{array} \right )| \ll{8vmhash}\ee
\[ q,r,s,t\in\complex\setminus\{0\}, \ p\in\complex\setminus\{1,-1\} \ \}  \]
where
\be Q=D(q)=diag(1,q,q,q^2) \ll{scalingm}\ee
This is the set of
the "eight--vertex" solutions up to the symmetries  
$R'=R^\#$.
The last two matrices can be transformed to $R_{3.1}(\frac{1-p}{1+p},\frac{1-p}{1+p},1)$ and $R_{1.4}(\frac{1+p}{1-p})$ by similarity \tfn s.

Below we display nontrivial invertible \sln s of \rf{xxzz}).
The list 
is complete up to the symmetries \rf{wxzsym})--\rf{dsym3}), \rf{condsym02}), \rf{condsym12}).
The ranges of parameters are
\[q,r,s,t,x\in\complex\setminus\{ 0 \}, \ \ 
\ u,v,a,b,c,d,e,f,g,h\in\complex, \ \ 
\ \epsilon,\epsilon_{1},\epsilon_{2}=\pm 1, \]
and ${\cal S},\ {\cal S}_{8V},\ {\cal S}_{5V},\ {\cal S}_{ST}$ are defined by \rf{ybesms})--\rf{stms}).
\vskip 1cm
{\bf List of nontrivial solutions of the equation \rf{xxzz}):}
\begin{itemize}
\item $X=X_{\ref{diag}}(a,b,c)$   
\\ The nontrivial solutions for general values of $a,b,c$ are
   \begin{enumerate}
   \item $Z\in {\cal S}_{5V}$
\vskip 2mm
For special values of $a,b,c$ we get the following solutions
   \item $(a=-1\wedge c=\pm b)$ or $(a=1\wedge c= -b)$ : $Z\in {\cal R}_{8V}$ 
   \item $a=1\wedge c= b$ : $Z\in {\cal S}\cup\{R_{1.2}^T(s)\}$.
 \end{enumerate}
\item $ X=X_{\ref{X.jj1}}(a,b,c) $   
\\ The nontrivial solutions for general values of $a,b,c$ are
   \begin{enumerate}
   \item $Z\in{\cal S}_{ST}$
   \end{enumerate}
There are no other solutions for special values of $a,b,c$.
\item $ X=X_{\ref{X.jj2}}(a,b,c) $   
\\ The nontrivial solutions for general values of $a,b,c$ are
   \begin{enumerate}
   \item $Z\in {\cal S}_{5V}$
\vskip 2mm
For special values of $a,b$ we get the following solutions
   \item $a=-1\wedge b=-c$ : $Z\in {\cal R}_{8V}$ 
   \item $a=1\wedge c= b$ : $ Z\in {\cal S}\cup\{R_{1.2}^T(s)\}$
 \end{enumerate}
\item $X=X_{\ref{X10}}(a,b,c)$   
\\ The nontrivial solutions for general values of $a,b,c$ are
   \begin{enumerate}
   \item $Z\in {\cal S}_{ST}$
   \end{enumerate}
There are no other solutions for special values of $a,b,c$.
\item $ X=X_{\ref{X8.3}}(a,b,c) $   
\\ There are no nontrivial solutions for general values of $a,b,c$. For special values we get
   \begin{enumerate}
   \item $b=-1, a=1/c\ :\ Z=R_{1.1}(i)$
   \item $b=-1, a=1/c\ :\ Z=R_{1.4}(x)\quad $
   \item $b=1\ :\ Z\in {\cal S}_{5V}$
   \item $b=1\wedge a=1/c\ :\ Z\in {\cal R}_{8V}$
   \end{enumerate}
\item $ X=X_{\ref{X8.4}}(a,b,c) $  
\\ There are no nontrivial solutions for general values of $a,b,c$. For special values we get
   \begin{enumerate}
   \item $a=b/c\ :\ Z=R_{1.1}(i)$
   \item $a=b/c\ :\ Z=R_{1.4}(x)\quad $
   \item $a=c/b\ :\ Z\in {\cal S}_{5V}$
   \item $a=-1,\  b=-c\ :\ Z\in {\cal R}_{8V}$
   \item $a=1,\  b=c\ :\ Z\in {\cal S}\cup\{R_{1.2}^T(s)\}$
   \end{enumerate}
\item $ X=X_{\ref{X8.5}}(a,b) $  
\\ There are no nontrivial solutions for general values of $a,b$. For special values we get
   \begin{enumerate}
   \item $a=b\ :\ Z\in {\cal S}_{ST}$
   \item $a=-b\ :\ Z=R_{3.1}(-1,-1,1)$
   \item $a=-b\ :\ Z=(S(x,y) \otimes S(x,y))^{-1} R_{0.2} (S(x,y) \otimes S(x,y))\quad $
\\ where $S(x,y)=\left ( \begin{array}{cc} x & 0 \\ y & 1 \end{array} \right ) $
   \end{enumerate}
\item $ X=X_{\ref{X8.7}}(a) $  
\\ The nontrivial solutions for general values of $a$ are
   \begin{enumerate}
     \item $Z= ( A \otimes A) R_{1.4}(x) ( A \otimes A)^{-1} \quad 
$
     \item $ Z= ( A \otimes A) R_{0.3} ( A \otimes A)^{-1} $
    where $
 A= \left ( \begin{array}{cc} \frac{1+\epsilon i}{\sqrt{2}} \sqrt{a} & \epsilon\sqrt{a} \\ \frac{1-\epsilon i}{\sqrt{2}} & i \end{array} \right )$
     \item $Z= ( B \otimes B) R_{3.1}(x,-x,1) ( B \otimes B)^{-1} \quad $
$, \ {\rm where} \
 B= \left ( \begin{array}{cc} -\sqrt{a}i & \sqrt{a}i \\ 1 & 1 \end{array} \right )$
   \end{enumerate}
There are no other solutions for special values of $a$.

 \item $ X=X_{\ref{h5}}(b,c,d) $  
\\ The nontrivial solutions for general values of $b,c,d$ are
   \begin{enumerate}
     \item $Z=R_{2.1}(c,b/d)$
\vskip 2mm
For special values of $b,c,d$ we get the following solutions ($D(q) $ is given by \rf{scalingm})).

     \item $b= c =1=-d : \quad 
     Z= D(q)R_{1.2}(x) D(q)^{-1}\quad$
     \item $b = d = -c=1  : \quad 
     Z =D(q)R_{1.2}^{at}(x)D(q)^{-1}\quad $
     \item $c = -d/b : \quad Z=R_{2.2}(c,x) \quad $
     \item $c=-1\wedge b = -d = \epsilon : \quad 
     Z=D(q) R_{0.2} D(q)^{-1}\quad $
     \item $b=1 \wedge c=d: \quad Z=(S(x) \otimes S(x)) R_{3.1}(d,1/d,1) (S(x) \otimes S(x))^{-1} \quad $
\\ where 
$ S(x)= \left ( \begin{array}{cc} {1-d} & 0 \\ x & x \end{array}
     \right ),\ d\neq 1 $
     \item $ b =  c =d=1  : \quad Z\in{\cal S}_{ST}$
     \item $c=1\wedge b=d = \epsilon : \quad Z=R_{0.1}$     
 \end{enumerate}    
 
\item $ X=X_{\ref{h6}}(b,c,d,g) $ 
\\ There are no nontrivial solutions for general values of parameters. For special values we get the following solutions ($D(q) $ is given by \rf{scalingm})).
   \begin{enumerate}
   \item $b=\epsilon d,\ c=\epsilon \ :\ Z=R_{3.1}(\epsilon,\epsilon,1)$
   \item $b=\epsilon d,\ c=-\epsilon \ :\ Z=R_{2.1}(-\epsilon,\epsilon)$
   \item $b=-d,\ c=1 \ :\quad Z =D(\sqrt{gx})R_{1.2}(-1+x(1-b^2))D(\sqrt{gx})^{-1}\quad $
   \item $b=d,\ c=-1 \ :\quad Z =D(\sqrt{gx})R_{1.2}^{at}(-1-x(1-b^2))D(\sqrt{gx})^{-1}\quad $
   \item $b=c=d=1,\ :\quad Z \in{\cal S}_{ST}$
   \item $b=1,\ c=d=-1\ :\quad Z =Z= (A(x,y) \otimes A(x,y))^{-1} R_{0.2} (A(x,y) \otimes A(x,y)) \quad $
\\where $ A(x,y)=\left ( \begin{array}{cc} x & 0 \\ y & 1 \end{array} \right ) $                
    \item $d=-b=-c=1\ :\quad Z =R_{0.2}$
    \item $c=-b=-d=1\ :\quad Z =R_{0.1}$
 \end{enumerate}

\item $ X=X_{\ref{X.30}}(a,b,c,d),\ c\neq 0$   
 \nogensol
   \item $a=c\ :  \quad Z =(A\otimes A)R_{3.1}(x,x,1)(A\otimes A)^{-1} \quad $ 
   \item $a=-c\ :  \quad Z =(A\otimes A)R_{1.4}(x)(A\otimes A)^{-1}\quad $
\\ where $A= \left ( \begin{array}{cc} \sqrt{b} & \sqrt{b} \\ -i & i \end{array}
     \right )  $
   \item 
   $a=c,\ d=b\ :  \quad Z = (B(y)\otimes B(y))R_{3.1}(x,x,1)(B(y)\otimes B(y))^{-1} \quad $
\\ where $B(y)= \left ( \begin{array}{cc} \sqrt{b}(1+\epsilon\sqrt{2})&y\sqrt{b}(1-\epsilon\sqrt{2})\\1&y\end{array}
     \right )  $
   \item 
   $a=c,\ d=b\ :  \quad Z = (C\otimes C)R_{1.4}(x)(C\otimes C)^{-1} \quad $ 
\\ where $C= \left ( \begin{array}{cc} \epsilon_1\sqrt{b}(\sqrt{2}+\epsilon_2)&\sqrt{b}\\ -\epsilon_1\epsilon_2&1+\epsilon_2\sqrt{2}\end{array}
     \right )  $
   \item 
   $a=c,\ d=b\ :  \quad Z = (D(x)\otimes D(x))R_{1.4}(1)(D(x)\otimes D(x))^{-1} \quad $ 
\\ where $D(x)= \left ( \begin{array}{cc} x\sqrt{b}&\sqrt{b}\\ \epsilon\sqrt{2}-x&\epsilon x\sqrt{2}-1\end{array}
     \right )  $
   \item 
   $a=c,\ d=b\ :  \quad Z = (E\otimes E)R_{1.4}(x)(E\otimes E)^{-1} \quad $ 
\\ where $E= \left ( \begin{array}{cc} i\sqrt{b}\epsilon_1&i\sqrt{ib}\epsilon_2\\1&-\epsilon_1\epsilon_2 \sqrt{i}\end{array}
     \right )  $
   \item 
   $a=c,\ d=b\ :  \quad Z = (F(y)\otimes F(y))R_{1.4}(1)(F(y)\otimes F(y))^{-1} \quad $
\\ where $F(y)= \left ( \begin{array}{cc} \sqrt{b}(\epsilon_1\sqrt{2}+y\epsilon_2)&\sqrt{b}\\y&-\epsilon_1 y\sqrt{2}- \epsilon_2\end{array}
     \right )  $
   \item 
   $a=-c,\ d=-b\ :  \quad Z = (G(y)\otimes G(y))R_{3.1}(x,-x,1)(G(y)\otimes G(y))^{-1} \quad $ 
\\ where $G(y)= \left ( \begin{array}{cc} \epsilon\sqrt{d}&-y\epsilon\sqrt{d}\\1&y\end{array}
     \right )  $
    \item   $a=-c,\ d=-b\ :  \quad Z = (H\otimes H)R_{0.3}(H\otimes H)^{-1} \quad $ 
\\ where $H= \left ( \begin{array}{cc} i\sqrt{b} & i\epsilon \sqrt{b}\\1&-\epsilon\end{array}     \right )  $
 \end{enumerate}     
 
\item $ X=X_{\ref{h7}}(a,b) $   
\\ The nontrivial solutions for general values of $a,b$ are
   \begin{enumerate}
     \item $Z=(A \otimes A)R_{2.2}(b/a,b/a)(A \otimes A)^{-1}$  
     \item $Z=(A \otimes A)R_{1.1}(a/b)(A \otimes A)^{-1}$  
    where $A=\left ( \begin{array}{cc} 1&1\\1&-1 \end{array}\right )  $
     \item $Z=( B \otimes B) R_{1.2}(b/a) ( B \otimes B)^{-1}$ 
     \\where $B=\left ( \begin{array}{cc} 1 & q\\1 & -q \end{array}  \right ),\ q^2={1-b/a}$
     \item $Z=( C \otimes C) R_{1.2}(-b/a) ( C \otimes C)^{-1}$ 
     \\where $
     C=\left ( \begin{array}{cc} 1& q\\-1&q\end{array}  \right ),\ q^2=1+b/a$
     \item $Z=( D \otimes D) R_{1.2}^T(-a/b) ( D \otimes D)^{-1}$ 
     \\where $
     D=\left ( \begin{array}{cc} q&1\\-q&1 \end{array}  \right ),\ q^2={-1-a/b}$
     \item $Z=( E \otimes E) R_{1.2}^T(a/b) ( E \otimes E)^{-1}$ 
     \\where $
     E=\left ( \begin{array}{cc} q&1\\q&-1 \end{array}  \right ),\ q^2={-1+a/b}$
     \item $Z=( F \otimes F) R_{0.3} ( F \otimes F)^{-1}$ where $
     F=\left ( \begin{array}{cc} \sqrt{i}&1\\ \sqrt{ib}&-\sqrt{b} \end{array}  \right )$
   \end{enumerate}
There are no other solutions for special values of $a,b$.

\item $ X=X_{\ref{octag}}(g) $ 
\\ There is no nontrivial solution of the equations \rf{xxzz}).
\item $ X=X_{\ref{X2.1}}(a,\epsilon,\frac{r+1}{r-1}) $ 
\\ The nontrivial solutions for general values of $a,\ r\neq 0,\pm 1$ are
   \begin{enumerate}
     \item $Z=(A \otimes A)R_{1.1}(\epsilon a/c)(A \otimes A)$  
    where $A=\left ( \begin{array}{cc} 1&1\\1&-1 \end{array}\right ),\ c^2=a^2-\frac{r-1}{r+1}  $
   \end{enumerate}
\item $ X=X_{\ref{X9.9}}(a) $ 
\\ The nontrivial solutions for general values of $a$ are
   \begin{enumerate}
     \item $Z=(A \otimes A)R_{3.1}(x,x,1)(A \otimes A)^{-1}$  
    \\ where $A=\left ( \begin{array}{cc} -q&q\\1&1 \end{array}\right ),\ q^2=a^2-1  $
\\ For special values of $a$ we get the following solutions 
    \item $a=\pm 1\ :\ Z=R_{0.1}$ 
    \item $a=0\ :\ Z=( B \otimes B) R_{3.1}(x,x,1) ( B \otimes B)^{-1}$ 
    \item $a=0\ :\ Z=( B \otimes B) R_{1.4}(x) ( B \otimes B)^{-1}$ 
     where $B=\left ( \begin{array}{cc} -1 & 1\\1 & 1 \end{array}  \right )$
    \item $a=0\ :\ Z=( C \otimes C) R_{1.4}(x) ( C \otimes C)^{-1}$ 
     where $C=\left ( \begin{array}{cc} -i & i\\1 & 1 \end{array}  \right )$
   \end{enumerate}

\item $ X=X_{\ref{h11.2}}(a,b,c,d) $ 
\nogensol  
     \item $a=d=0 : \quad Z\in{\cal S}_{5V}$
    \item $c=b(1+a-d) : \quad Z=R_{1.3}(d-a-1)$
 \end{enumerate}    

\item $ X=X_{\ref{h11.1}}(a,b) $ 
\\ There is no nontrivial solution of the equations \rf{xxzz}).
\item $ X=X_{\ref{X9.4}}(b) $ 
\\ The nontrivial solutions for  $b\neq 0,\epsilon$ are
   \begin{enumerate}
     \item $Z= 
(A \otimes A) R_{1.4}(x) ( A \otimes A)^{-1}, \quad $
\\     where $A=\left ( \begin{array}{cc} i & -ip\\q&pq \end{array}  \right ),\ q^2=1-b^2,\ p^4=(b-iq)^2.$
\\For special values of $b$ we get 
\item $b=\epsilon\ :  Z=( B(x) \otimes B(x)) R_{0.2} ( B(x) \otimes B(x))^{-1}$ 
     \\ where $B(x)=\left ( \begin{array}{cc} -1 & x \\ 4\epsilon & 0 \end{array}  \right )$ 
\item $b=0$ see $X_{\ref{X8.7}}(a=1)$
 \end{enumerate}    

\item $ X=X_{\ref{X10.5}}(a,b,c,\epsilon) $ 
\nogensol
     \item $\epsilon=1,b=0 : \quad Z\in {\cal S}_{5V}$
     \item $\epsilon=-1,\ a=\pm 1 : \quad Z=R_{1.1}(i)$
     \item $\epsilon=-1,\ a=\pm 1 : \quad Z=R_{1.4}(x)$
     \item $\epsilon=-1,\ c=ab : \quad Z\in {\cal S}_{5V}$
     \item $\epsilon=-1,\ a=-1,\ c=-b : \quad Z\in {\cal R}_{8V}$
     \item $\epsilon=-1,\ a=1,\ c=b : \quad Z\in {\cal S}\cup\{R_{1.2}^T(s)\}$
 \end{enumerate}    

\item $ X=X_{\ref{X10.6}}(a,b,c) $ 
\nogensol
  \item $ b = 0 : \quad     Z= R_{3.1}(-1,-1,1)$
     \item $ b = 0 : \quad     Z= (A(x,y) \otimes A(x,y))^{-1} R_{0.2} (A(x,y) \otimes A(x,y)) \quad $
\\where $ A(x,y)=\left ( \begin{array}{cc} x & 0 \\ y & 1 \end{array} \right ) $                
     \item $ a = c b : \quad Z\in {\cal S}_{ST}$
 \end{enumerate}    

\item $ X=X_{\ref{X1.4}}(a,b) $ 
 \nogensol
   \item $ b = 0 : \quad Z= (S \otimes S) R_{0.3} (S \otimes S)^{-1} \quad {\rm where} \; S=\left ( \begin{array}{cc}1& \sqrt{i} \\ 1 & -\sqrt{i}  \end{array} \right )  $                
  \item $ a = b =0 : \quad Z= (S \otimes S) R_{3.1}(x,-x,1) (S \otimes S)^{-1} \quad  \; {\rm where} \; S=\left ( \begin{array}{cc}1& 1 \\ -1 & 1  \end{array} \right ) $                
  \item $ a = b =0 : \quad Z= (S \otimes S) R_{1.4}(x) (S \otimes S)^{-1} \quad  \; {\rm where} \; S=\left ( \begin{array}{cc}1& \sqrt{i} \\ -1 & \sqrt{i}  \end{array} \right ) $                
  \end{enumerate}
\item $ X=X_{\ref{X8.15}}(a,b,c,d,e,f,g,h) $ and $ X=X_{\ref{X8.16}}(a,b,c,d,e,f,g,h) $.
\\ Solving the system \rf{xxzz}) for these two matrices is rather difficult. We can, however, simplify the task by the following way: 

Both the matrices solve the equation \rf{wxx}) only for $W=\unit$, so that any $Z$ that solve \rf{xxzz}) for $X_{\ref{X8.15}},X_{\ref{X8.16}}$ gives solution $(\unit,X,Z)$ of the WXZ system. This solution is equivalent to $(Z^\pm,X^+,\unit)$ due to \rf{dsym3}). It is always possible to transform $Z^\pm$ by WXZ symmetries to the form contained in ${\cal S}$. It means that for $Z\neq \unit$ we can find the solution of the WXZ by symmetries investigating $W\neq \unit$. It remains to investigate the case $Z=\unit$ or, in other words, to solve the equation $[X,X,\unit]=0$ where $X$ is of the form $X_{\ref{X8.15}}$ or $X_{\ref{X8.16}}$. It turns out that the resulting matrices can be always transformed by the symmetry \tfn {} $X'=(T\otimes\unit) X (T\otimes\unit)^{-1}$ to $X_1,X_2,X_3,X_4.$
\end{itemize}

\subsection{List of solutions of the WXZ system}
Combining the results of the subsections \ref{wxxsolns} and \ref{xxzzsolns} we get the list of solutions of the WXZ system  \rf{www})--\rf{xxz}) presented below. By construction the list is complete up to the symmetries \rf{wxzsym})--\rf{dsym3}) and \rf{condsym01}),\rf{condsym1}),\rf{condsym2}).

Let us remind that $a,b,c,d,g,h\in\complex$, $\epsilon,\epsilon_1,\epsilon_2,\epsilon_3=\pm 1$, and ${\cal S},\ {\cal S}_{8V},\ {\cal S}_{5V},\ {\cal S}_{ST}$ are defined by \rf{ybesms})--\rf{stms}).
\vskip 1cm{\bf Generic solutions:}
\begin {enumerate}
 \item $W\in{\cal S},\ \ X=\unit,\ \ Z\in {\cal S}\cup\{R_{1.2}^T(s)\}$
 \item $W\in{\cal S}_{8V},\ \ X=diag(1,1,1,-1),\ \ Z\in {\cal S}_{8V}\cup\{R_{1.2}^T(s)\}$
\item $W\in{\cal S}_{8V},\ \ X=\left ( \begin{array}{cccc}
         1 & 0 & 0 & 0 \\
         0 & 1 & 0 & 0 \\
         0 & 0 & 0 & 1 \\
         0 & 0 & 1 & 0
       \end{array} \right ) ,\ \ Z=R_{3.1}(x,x,1)$
\item $W\in{\cal S}_{8V},\ \ X=\left ( \begin{array}{cccc}
         1 & 0 & 0 & 0 \\
         0 & 1 & 0 & 0 \\
         0 & 0 & 0 & 1 \\
         0 & 0 & 1 & 0
       \end{array} \right ) ,\ \ Z=R_{1.4}(x)$
\item $W\in{\cal S}_{5V},\ \ X=diag(1,1,1,a),\ \ Z\in{\cal S}_{5V}\ {\rm or}\ \ Z=P$
\item $W\in{\cal S}_{5V},\ \ X=\left ( \begin{array}{cccc}
         1 & 0 & 0 & 0 \\
         0 & 1 & 0 & 0 \\
         0 & 0 & 1 & 0 \\
         0 & 0 & a & 1
       \end{array} \right ) ,\ \ Z\in{\cal S}_{ST}$ or $Z=P$
       
\item $W\in{\cal S}_{ST},\ \ X=\left ( \begin{array}{cccc}
  1 & 0 & 0 & 0 \\
  0 & 1 & 0 & 0 \\
  0 & 0 & 1 & 0 \\
  a & 0 & 0 & 1
  \end{array} \right ) ,\ \ Z\in{\cal S}_{ST}$ or $Z=P$
\vskip 1cm
{\bf Non--generic solutions:}
\item $W=R_{3.1}(s,s,1),\ \ X=\left ( \begin{array}{cccc}
  1 & 0 & 0 & 0 \\
  0 & 0 & 0 & 1 \\
  0 & 0 & b & 0 \\
  0 & 1 & 0 & 0
       \end{array} \right ) ,\ 
  \ Z=P$
\item $W=R_{3.1}(s,s,1),\ \ X=\left ( \begin{array}{cccc}
  1 & 1 & 1 & -1 \\
  1 & 1 & -1 & 1 \\
  1 & -1 & 1 & 1 \\
  -1 & 1 & 1 & 1
       \end{array} \right ) ,\ 
  \ Z=R_{3.1}(x,x,1)$
\item $W=R_{3.1}(s,s,1),\ \ X=\left ( \begin{array}{cccc}
  1 & 1 & 1 & -1 \\
  1 & 1 & -1 & 1 \\
  1 & -1 & 1 & 1 \\
  -1 & 1 & 1 & 1
       \end{array} \right ) ,\ 
  \ Z=R_{1.4}(x)$
\item $W=R_{3.1}(s,s,1),\ \ X=\left ( \begin{array}{cccc}
  1 & 0 & 0 & 0 \\
  0 & 0 & 0 & 1 \\
  0 & 0 & -1 & 0 \\
  0 & 1 & 0 & 0
       \end{array} \right ) ,\ 
  \ Z=R_{1.1}(i)$
\item $W=R_{3.1}(s,s,1),\ \ X=\left ( \begin{array}{cccc}
  1 & 0 & 0 & 0 \\
  0 & 0 & 0 & 1 \\
  0 & 0 & -1 & 0 \\
  0 & 1 & 0 & 0
       \end{array} \right ) ,\ 
  \ Z=R_{1.4}(x)$
\item $W=R_{3.1}(s,s,1),\ \ X=\left ( \begin{array}{cccc}
  0 & 0 & 1 & 0 \\
  0 & 0 & 0 & a \\
  1 & 0 & 0 & 0 \\
  0 & 1 & 0 & 0
       \end{array} \right ) ,\ 
  \ Z=P$
\item $W=R_{3.1}(s,s,1),\ \ X=\left ( \begin{array}{cccc}
  0 & 0 & 1 & 0 \\
  0 & 0 & 0 & a \\
  1 & 0 & 0 & 0 \\
  0 & 1/a & 0 & 0
       \end{array} \right ) ,\ 
  \ Z=R_{1.1}(i)$
\item $W=R_{3.1}(s,s,1),\ \ X=\left ( \begin{array}{cccc}
  0 & 0 & 1 & 0 \\
  0 & 0 & 0 & a \\
  1 & 0 & 0 & 0 \\
  0 & 1/a & 0 & 0
       \end{array} \right ) ,\ 
  \ Z=R_{1.4}(x)$
\item $W=R_{3.1}(s,s,1),\ \ X=\left ( \begin{array}{cccc}
  0 & 0 & 1 & 0 \\
  0 & 0 & 1 & 1 \\
  1 & 0 & 0 & 0 \\
  c & 1 & 0 & 0
       \end{array} \right ) ,\ 
  \ Z=P$
\item $W=R_{3.1}(s,s,1),\ \ X=\left ( \begin{array}{cccc}
  0 & 0 & 1 & 0 \\
  0 & 0 & 1 & 1 \\
  1 & 0 & 0 & 0 \\
  -1 & 1 & 0 & 0
       \end{array} \right ) ,\ 
  \ Z=R_{3.1}(-1,-1,1)$
\item $W=R_{3.1}(s,s,1),\ \ X=\left ( \begin{array}{cccc}
  0 & 0 & x & 0 \\
  0 & 0 & 1 & x \\
  x & 0 & 0 & 0 \\
  -1 & x & 0 & 0
       \end{array} \right ) ,\ 
  \ Z=R_{0.2}$
\item $W=R_{3.1}(s,-s,1),\ \ X=\left ( \begin{array}{cccc}
  0 & 0 & 1 & 0 \\
  0 & 0 & 0 & -1 \\
  0 & a & 0 & 0 \\
  1 & 0 & 0 & 0
       \end{array} \right ) ,\ 
  \ Z=P$
\item $W=R_{3.1}(s,-s,1),\ \ X=\left ( \begin{array}{cccc}
  0 & 0 & 0 & 1 \\
  0 & 0 & 1 & 0 \\
  0 & q & 0 & 0 \\
  -q & 0 & 0 & 0
       \end{array} \right ) ,\ 
  \ Z=R_{3.1}(x,-x,1)$
\item $W=R_{3.1}(s,-s,1),\ \ X=\left ( \begin{array}{cccc}
  0 & 0 & 0 & 1 \\
  0 & 0 & i\epsilon & 0 \\
  0 & i\epsilon  q& 0 & 0 \\
  q & 0 & 0 & 0
       \end{array} \right ) ,\ 
  \ Z=R_{1.4}(x)$
\item $W=R_{3.1}(s,-s,1),\ \ X=\left ( \begin{array}{cccc}
  0 & 0 & 0 & 1 \\
  0 & 0 & i\epsilon & 0 \\
  0 & i\epsilon q & 0 & 0 \\
  q & 0 & 0 & 0
       \end{array} \right ) ,\ 
  \ Z=R_{0.3}$

\item $W=R_{3.1}(\epsilon_1,\epsilon_1,1),\ \ 
X=\left ( \begin{array}{cccc}
  1 & 0 & 0 & 0 \\
  0 & \epsilon_1 & 1 & 0 \\
  0 & 0 & \epsilon_2 & 0 \\
  g & 0 & 0 & \epsilon_1\epsilon_2
       \end{array} \right ) ,\ 
  \ Z=R_{3.1}(\epsilon_2,\epsilon_2,1)$
\item $W=R_{3.1}(\epsilon,\epsilon,1),\ \ 
X=\left ( \begin{array}{cccc}
  1 & 0 & 0 & 0 \\
  0 & \epsilon & 1 & 0 \\
  0 & 0 & 1 & 0 \\
  g & 0 & 0 & \epsilon
       \end{array} \right ) ,\ 
  \ Z=R_{0.1}$
\item $W=R_{3.1}(\epsilon,\epsilon,1),\ \ 
X=\left ( \begin{array}{cccc}
  1 & 0 & 0 & 0 \\
  0 & \epsilon & 1 & 0 \\
  0 & 0 & -1 & 0 \\
  g & 0 & 0 & -\epsilon
       \end{array} \right ) ,\ 
  \ Z=R_{0.2}$
\item $W=R_{3.1}(-1,-1,1),\ \ 
X=\left ( \begin{array}{cccc}
  0 & ab & c & 0 \\
  a & 0 & 0 & -c \\
  d & 0 & 0 & b \\
  0 &-d & 1 & 0
       \end{array} \right ) ,\ 
  \ Z=P$
\item $W=R_{3.1}(1,1,1),\ \ X=\left ( \begin{array}{cccc}
         1 & 0 & 0 & 0 \\
         a & 1 & b & 0 \\
         c & 0 & d & 0 \\
         g & c & h & d
       \end{array} \right ) ,\ 
  \ Z=P$
\item $W=R_{3.1}(1,1,1),\ \ X=\left ( \begin{array}{cccc}
         1 & 0 & 0 & 0 \\
         0 & a & 0 & b \\
         c & 0 & d & 0 \\
         0 & g & 0 & h
       \end{array} \right ) ,\ 
  \ Z=P$

\item $W=R_{2.1}(r,s),\ \ X=\left ( \begin{array}{cccc}
  1 & 0 & 0 & 0 \\
  0 & r^{-1} & 0 & 0 \\
  0 & 0 & a & 0 \\
  1 & 0 & 0 & sa
       \end{array} \right ) ,\ 
  \ Z=P$
\item $W=R_{2.1}(r,s),\ \ X=\left ( \begin{array}{cccc}
  1 & 0 & 0 & 0 \\
  0 & r^{-1} & 0 & 0 \\
  0 & 0 & a & 0 \\
  1 & 0 & 0 & sa
       \end{array} \right ) ,\ 
  \ Z=R_{2.1}(a,(rsa)^{-1})$
\item $W=R_{2.1}(-1/s,s),\ \ X=\left ( \begin{array}{cccc}
  1 & 0 & 0 & 0 \\
  0 & -s & 0 & 0 \\
  0 & 0 & a & 0 \\
  1 & 0 & 0 & sa
       \end{array} \right ) ,\ 
  \ Z=R_{2.2}(a,x)$
\item $W=R_{2.1}(\epsilon,-\epsilon),\ \ X=\left ( \begin{array}{cccc}
  1 & 0 & 0 & 0 \\
  0 & \epsilon & 0 & 0 \\
  0 & 0 & 1 & 0 \\
  q & 0 & 0 & -\epsilon
       \end{array} \right ) ,\ 
  \ Z=R_{1.2}(x)$
\item $W=R_{2.1}(s,s),\ \ X=\left ( \begin{array}{cccc}
  1 & 0 & 0 & 0 \\
  0 & s^{-1} & 1 & 0 \\
  0 & 0 & a & 0 \\
  g & 0 & 0 & sa
       \end{array} \right ) ,\ 
  \ Z=P$
\item $W=R_{2.1}(i\epsilon_1,i\epsilon_1),\ \ X=\left ( \begin{array}{cccc}
  1 & 0 & 0 & 0 \\
  0 & -i\epsilon_1 & 1 & 0 \\
  0 & 0 & \epsilon_2 & 0 \\
  g & 0 & 0 & i\epsilon_1\epsilon_2
       \end{array} \right ) ,\ 
  \ Z=R_{2.1}(\epsilon_2,-\epsilon_2)$
\item $W=R_{2.1}(i\epsilon,i\epsilon),\ \ X=\left ( \begin{array}{cccc}
  1 & 0 & 0 & 0 \\
  0 & -i\epsilon & q & 0 \\
  0 & 0 & 1& 0 \\
  a & 0 & 0 & i\epsilon
       \end{array} \right ) ,\ 
  \ Z=R_{1.2}(-1+\frac{2}{aq})$
\item $W=R_{2.1}(1,-1),\ \ X=\left ( \begin{array}{cccc}
  a & 0 & 0 & 0 \\
  0 & a & 0 & 0 \\
  1 & 0 & 0 & b \\
  0 &-1 & b & 0
       \end{array} \right ) ,\ 
  \ Z=P$
\item $W=R_{2.1}(1,-1),\ \ X=\left ( \begin{array}{cccc}
  a & 0 & 0 & 0 \\
  0 & a & 0 & 0 \\
  0 & 1 & b & 0 \\
  1 & 0 & 0 & -b
       \end{array} \right ) ,\ 
  \ Z=R_{2.2}(b/a,b/a)$
\item $W=R_{2.1}(1,-1),\ \ X=\left ( \begin{array}{cccc}
  a & 0 & 0 & 0 \\
  0 & a & 0 & 0 \\
  0 & 1 & b & 0 \\
  1 & 0 & 0 & -b
       \end{array} \right ) ,\ 
  \ Z=R_{1.1}(a/b)$
\item $W=R_{2.1}(\epsilon,-\epsilon),\ \ X=\left ( \begin{array}{cccc}
  a & 0 & 0 & 0 \\
  0 & \epsilon a & 0 & 0 \\
  0 & \sqrt{\frac{a-b}{\epsilon a}} & b & 0 \\
  \sqrt{\frac{\epsilon a}{a-b}} & 0 & 0 & -\epsilon b
       \end{array} \right ) ,\ 
  \ Z=R_{1.2}(b/a)$
\item $W=R_{2.1}(1,-1),\ \ X=\left ( \begin{array}{cccc}
  a & 0 & 0 & 0 \\
  0 & a & 0 & 0 \\
  0 & 1/\sqrt{i} & b & 0 \\
  \sqrt{i} & 0 & 0 & -b
       \end{array} \right ) ,\ 
  \ Z=R_{0.3}$

\item $W=R_{2.2}(r,s),\ \ X=\left ( \begin{array}{cccc}
  1 & 0 & 0 & 0 \\
  0 & r^{-1} & 0 & 0 \\
  0 & 0 & a & 0 \\
  1 & 0 & 0 & -ar^{-1}
       \end{array} \right ) ,\ 
  \ Z=P$
\item $W=R_{2.2}(r,s),\ \ X=\left ( \begin{array}{cccc}
  1 & 0 & 0 & 0 \\
  0 & r^{-1} & 0 & 0 \\
  0 & 0 & a & 0 \\
  1 & 0 & 0 & -ar^{-1}
       \end{array} \right ) ,\ 
  \ Z=R_{2.2}(a,x)$
\item $W=R_{2.2}(\epsilon,s),\ \ X=\left ( \begin{array}{cccc}
  1 & 0 & 0 & 0 \\
  0 & \epsilon & 0 & 0 \\
  0 & 0 & 1 & 0 \\
  q & 0 & 0 & -\epsilon
       \end{array} \right ) ,\ 
  \ Z=R_{1.2}(x)$
  \item $W=R_{2.2}(s,s),\ \ X=\left ( \begin{array}{cccc}
  1 & 0 &0 & 0 \\
  0 & s^{-1} & 1 & 0 \\
  0 & 0 & a & 0 \\
  g & 0 & 0 & -a/s
       \end{array} \right ) ,\ 
  \ Z=P$
  \item $W=R_{2.2}(s,s),\ \ X=\left ( \begin{array}{cccc}
  1 & 0 &0 & 0 \\
  0 & s^{-1} & q & 0 \\
  0 & 0 & 1 & 0 \\
  g & 0 & 0 & -1/s
       \end{array} \right ) ,\ 
  \ Z=R_{1.2}(\frac{s^2-1}{qgs^2}-1)$

\item $W=R_{1.1}(r),\ \ X=\left ( \begin{array}{cccc}
  1 & 0 & 0 & 0 \\
  0 & r & 1 & 0 \\
  0 & 0 & -\epsilon & 0 \\
  \epsilon & 0 & 0 & \epsilon r
       \end{array} \right ) ,\ 
  \ Z=P$
\item $W=R_{1.1}(r),\ \ X=\left ( \begin{array}{cccc}
  1 & 0 & 0 & 0 \\
  0 & r & q & 0 \\
  0 & 0 & 1 & 0 \\
  -q & 0 & 0 & -r
       \end{array} \right ) ,\ 
  \ Z=R_{1.2}(\frac{r^2-1}{q^2}-1)$
\item $W=R_{1.1}(r),\ \ X=\left ( \begin{array}{cccc}
\epsilon a\frac{r+1}{r-1} & \epsilon c & 1 & 0 \\
     \epsilon c& \epsilon a \frac{r+1}{r-1} & 0 & -1 \\
         \epsilon & 0 & a & c\frac{r+1}{r-1}  \\
         0 & -\epsilon & c\frac{r+1}{r-1} & a
       \end{array} \right ) ,\ c^2=a^2-\frac{r-1}{r+1}
      ,\ 
  \\ Z=P$
\item $W=R_{1.1}(r),\ \ X=\left ( \begin{array}{cccc}
   \epsilon (a\frac{r+1}{r-1} +c) & 0 & 0 & 1 \\
   0 & \epsilon (a\frac{r+1}{r-1}-c) & 1 & 0  \\
   0 &      \epsilon &(a+c\frac{r+1}{r-1} )& 0  \\
   \epsilon & 0 & 0 & (a-c\frac{r+1}{r-1})
       \end{array} \right ) ,
        \\ c^2=a^2-\frac{r-1}{r+1},\ Z=R_{1.1}(\epsilon a/c),$ 
\item $W=R_{1.1}(i),\ \ X=\left ( \begin{array}{cccc}
  a & 1-a^2 & 0 & 0 \\
  1 & -a & 0 & 0 \\
  0 & 0 & 1 & 0 \\
  0 & 0 & 0 & -1
       \end{array} \right ) ,\ 
  \ Z=P$
\item $W=R_{1.1}(i),\ \ X=\left ( \begin{array}{cccc}
  \epsilon & 0 & 0 & 0 \\
  1 & -\epsilon & 0 & 0 \\
  0 & 0 & 1 & 0 \\
  0 & 0 & 0 & -1
       \end{array} \right ) ,\ 
  \ Z=R_{0.1}$
\item $W=R_{1.1}(i),\ \ X=\left ( \begin{array}{cccc}
  1 & 0 & 0 & 0 \\
  0 & -1 & 0 & 0 \\
  0 & 0 & 0 & 1 \\
  0 & 0 & 1 & 0
       \end{array} \right ) ,\ 
  \ Z=R_{1.4}(x)$
\item $W=R_{1.1}(i),\ \ X=\left ( \begin{array}{cccc}
  0 & -i & 0 & 0 \\
  i & 0 & 0 & 0 \\
  0 & 0 & 0 & 1 \\
  0 & 0 & 1 & 0
       \end{array} \right ) ,\ 
  \ Z=R_{1.4}(x)$
\item $W=R_{1.2}(s),\ \ X=\left ( \begin{array}{cccc}
  1 & 0 & 0 & 0 \\
  0 & 1 & 0 & 0 \\
  0 & 0 & \epsilon & 0 \\
  1 & 0 & 0 & -\epsilon
  \end{array} \right ) ,\   \ Z=P$
\item $W=R_{1.2}(s),\ \ X=\left ( \begin{array}{cccc}
  1 & 0 & 0 & 0 \\
  0 & 1 & 0 & 0 \\
  0 & 0 & 1 & 0 \\
  q & 0 & 0 & -1
       \end{array} \right ) ,\ 
  \ Z=R_{1.2}(x)$
\item $W=R_{1.2}(s),\ \ X=\left ( \begin{array}{cccc}
  1 & 0 & 0 & 0 \\
  0 & 1 & 0 & 0 \\
  0 & q & 1 & 0 \\
  0 & 0 & 0 & -1
       \end{array} \right ) ,\ 
  \ Z=R_{1.2}^T(x)$
\item $W=R_{1.2}(s),\ \ X=\left ( \begin{array}{cccc}
  1 & 0 & 0 & 0 \\
  0 & s^{-1} & s-1 & 0 \\
  0 & 0 & \epsilon & 0 \\
  {\epsilon} & 0 & 0 & -\epsilon s^{-1}
       \end{array} \right ) ,\ 
  \ Z=P$
\item $W=R_{1.2}(s),\ \ X=\left ( \begin{array}{cccc}
  1 & 0 & 0 & 0 \\
  0 & s^{-1} & a(s-1) & 0 \\
  0 & 0 & 1 & 0 \\
  a & 0 & 0 & -s^{-1}
       \end{array} \right ) ,\ 
  \ Z=R_{1.2}(\frac{s+1}{a^2s^2}-1)$
\item $W=R_{1.2}(s),\ \ X=\left ( \begin{array}{cccc}
  s^{-1} & 0 & 0 & a(1-s)  \\
  0 & 1 & 0 & 0 \\
  0 & a & s^{-1} & 0 \\
  0 & 0 & 0 & -1 
       \end{array} \right ) ,\ 
  \ Z=R_{1.2}^T(\frac{s+1}{a^2s^2}-1)$
\item $W=R_{1.2}(s),\ \ X=\left ( \begin{array}{cccc}
  a & 0 & 0 & 0 \\
  0 & a & 0 & 0 \\
  1 & 0 & 0 & 1 \\
  0 &-1 & a^2 -s-1 & 0
       \end{array} \right ) ,\ 
  \ Z=P$
\item $W=R_{1.2}(s),\ \ X=\left ( \begin{array}{cccc}
  a & 0 & 0 & 0 \\
  0 & a & 0 & 0 \\
  0 & 1/\sqrt{i} & \sqrt{a^2-s-1} & 0 \\
  \sqrt{i}&0 & 0 & -\sqrt{a^2-s-1}
       \end{array} \right ) ,\ 
  \\ Z=R_{0.3}$
\item $W=R_{1.3}(u),\ \ X=\left ( \begin{array}{cccc}
         1 & 0 & 0 & 0 \\
         a & 1 & 0 & 0 \\
         b & 0 & 1 & 0 \\
         0 & b-u-1 & -ua & 1
       \end{array} \right ) ,\ 
  \ Z=P$
\item $W=R_{1.3}(u),\ \ X=\left ( \begin{array}{cccc}
         1 & 0 & 0 & 0 \\
         -1/b & 1 & 0 & 0 \\
         b & 0 & 1 & 0 \\
         0 & b-u-1 & u/b & 1
       \end{array} \right ) ,\ 
  \ Z=R_{1.3}(\frac{u+1}{b}-1)$
\item $W=R_{1.3}(1),\ \ X=\left ( \begin{array}{cccc}
         1 & 0 & 0 & 0 \\
         a & 1 & 1 & 0 \\
         a+b+1 & 0 & 1 & 0 \\
         0 & a+b-1 & b & 1
       \end{array} \right ) ,\ 
  \ Z=P$

\item $W=R_{1.4}(t),\ \ X=\left ( \begin{array}{cccc}
  1 & 0 & 0 & 0 \\
  0 & 0 & 0 & a \\
  0 & 0 & \epsilon & 0 \\
  0 & a & 0 & 0
       \end{array} \right ) ,\ 
  \ Z=P$
\item $W=R_{1.4}(t),\ \ X=\left ( \begin{array}{cccc}
  0 & 0 & 1 & 0 \\
  0 & 0 & 0 & a \\
  1 & 0 & 0 & 0 \\
  0 & -a & 0 & 0
       \end{array} \right ) ,\ 
  \ Z=P$
\item $W=R_{1.4}(t),\ \ X=\left ( \begin{array}{cccc}
  0 & 0 & 1 & 0 \\
  0 & 0 & 0 & i \\
  1 & 0 & 0 & 0 \\
  0 & -i & 0 & 0
       \end{array} \right ) ,\ 
  \ Z=R_{1.4}(x)$
\item $W=R_{1.4}(t),\ \ X=\left ( \begin{array}{cccc}
  a & 1-a^2 & 0 & 0 \\
  1 & -a & 0 & 0 \\
  0 & 0 & 1 & 0 \\
  0 & 0 & 0 & -1
       \end{array} \right ) ,\ 
  \ Z=P$
\item $W=R_{1.4}(t),\ \ X=\left ( \begin{array}{cccc}
  \epsilon & 0 & 0 & 0 \\
  1 & -\epsilon & 0 & 0 \\
  0 & 0 & 1 & 0 \\
  0 & 0 & 0 & -1
       \end{array} \right ) ,\ 
  \ Z=R_{0.1}$
\item $W=R_{1.4}(t),\ \ X=\left ( \begin{array}{cccc}
  1 & 0 & 0 & 0 \\
  0 & -1 & 0 & 0 \\
  0 & 0 & 0 & 1 \\
  0 & 0 & 1 & 0 \\
       \end{array} \right ) ,\ 
  \ Z=R_{1.4}(x)$
\item $W=R_{1.4}(t),\ \ X=\left ( \begin{array}{cccc}
  1 & 0 & 0 & 0 \\
  0 & -1 & 0 & 0 \\
  0 & 0 & 0 & i \\
  0 & 0 & -i & 0 \\
       \end{array} \right ) ,\ 
  \ Z=R_{1.4}(x)$
\item $W=R_{1.4}(t),\ \ X=\left ( \begin{array}{cccc}
  0 & 0 & 1 & 0 \\
  0 & 0 & 0 & -1 \\
  a & 1 & 0 & 0 \\
  1-a^2 & -a & 0 & 0
       \end{array} \right ) ,\ 
  \ Z=P$
\item $W=R_{1.4}(t),\ \ X=\left ( \begin{array}{cccc}
  0 & 0 & 0 & p \\
  0 & 0 & 1/p &  0\\
  0 &  \epsilon/p & 0 & 0\\
 \epsilon p &  0 & 0 & 0
       \end{array} \right ) ,
  \ Z=R_{1.4}(x)$
\item $W=R_{1.4}(t),\ \ X=\left ( \begin{array}{cccc}
  0 & 0 & 1 & 0 \\
  0 & 0 & x & -1 \\
  1 & 0 & 0 & 0 \\
  -x & -1& 0 & 0
       \end{array} \right ) ,\ 
  \ Z=R_{0.2}$

\item $W=R_{1.4}(x),\ \ X=\left ( \begin{array}{cccc}
  0 & 0 & 0 & 1 \\
  0 & 0 & i\epsilon & 0 \\
  0 & i\epsilon & 0 & 0 \\
  1 & 0 & 0 & 0
       \end{array} \right ) ,\ 
  \ Z=R_{1.4}(x)$
\item $W=R_{1.4}(x),\ \ X=\left ( \begin{array}{cccc}
  0 & 0 & 0 & 1 \\
  0 & 0 & i\epsilon & 0 \\
  0 & i\epsilon  & 0 & 0 \\
  1 & 0 & 0 & 0
       \end{array} \right ) ,\ 
  \ Z=R_{0.3}$
\item $W=R_{0.1},\ \ X=\left ( \begin{array}{cccc}
  1 & 0 & 0 & 0 \\
  0 & a & 0 & 0 \\
  b & 0 & \epsilon & 0 \\
  0 & c & 0 & -a
       \end{array} \right ) ,\ 
  \ Z=P$
\item $W=R_{0.1},\ \ X=\left ( \begin{array}{cccc}
  1 & 0 & 0 & 0 \\
  b & 1 & 0 & 0 \\
  c & 0 & -1 & 0 \\
  a & c & -b & -1
       \end{array} \right ) ,\ 
  \ Z=P$
\item $W=R_{0.1},\ \ X=\left ( \begin{array}{cccc}
  1 & 0 & 0 & 0 \\
  0 & 1 & 0 & 0 \\
  c & 0 & -1 & 0 \\
  a & c & 0 & -1
       \end{array} \right ) ,\ 
  \ Z=R_{0.2}$
\item $W=R_{0.2},\ \ X=\left ( \begin{array}{cccc}
  1 & 0 & 0 & 0 \\
  0 & -1 & 0 & 0 \\
  0 & 0 & \epsilon & 0 \\
  1 & 0 & 0 & -\epsilon
       \end{array} \right ) ,\ 
  \ Z=P$

\item $W=R_{0.2},\ \ X=\left ( \begin{array}{cccc}
  1 & 0 & 0 & 0 \\
  0 & -1 & 0 & 0 \\
  0 & 0 & -1 & 0 \\
  q & 0 & 0 & 1
       \end{array} \right ) ,\ 
  \ Z=R_{0.2}$
\item $W=R_{0.2},\ \ X=\left ( \begin{array}{cccc}
  0 & \epsilon b & 0 & 0 \\
  \epsilon & 0 & 0 & 0 \\
  1 & 0 & 0 & b \\
  0 &-1 & 1 & 0
       \end{array} \right ) ,\ 
  \ Z=P$

\item $W=R_{0.3},\ \ X=\left ( \begin{array}{cccc}
  1 & 0 & 0 & 0 \\
  0 & a & 1 & 0 \\
  0 & 0 & -1 & 0 \\
  i & 0 & 0 & a
       \end{array} \right ) ,\ 
  \ Z=P$
\item $W=R_{0.3},\ \ X=\left ( \begin{array}{cccc}
  a & b & 1 & 0 \\
  b & a & 0 & -1 \\
  0 & i & -b & -a \\
  -i & 0 & -a & -b
       \end{array} \right ) ,\ 
  \ Z=P$
\item $W=R_{0.3},\ \ X=\left ( \begin{array}{cccc}
  a & 0 & 0 & i \\
  0 & a & 1 & 0 \\
  0 & 1 & -a & 0 \\
  i & 0 & 0 & a
       \end{array} \right ) ,\ 
  \ Z=R_{0.3}$

\item $W=R_{0.4}=P,\ \ X={\rm arbitrary\ matrix},\ \ Z=P$

\end{enumerate}

\section {Conclusion}
It is possible to solve completely the system of equations \rf{www})--\rf{xxz}) in the dimension two with the assistence of the computer programs for symbolic manipulations. Very important tool for both solving and classification of solutions is the use of symmetries.

The invertible solutions of the WXZ system are classified in this paper  and the authors  believe that that they were careful enough so that the solution is complete up to the symmetries of the system \rf{wxzsym})--\rf{condsym2}). The rigorous proof of the completeness of the solution set is rather difficult task. It is not known  for pure \YBE {} either to the best knowledge of the authors.

The number of the solutions is rather large so that it is not possible to investigate in detail all the quantum doubles generated by the Vladimirov's procedure. It is therefore necessary to formulate supplementary conditions for the solutions that guarantee some suitable properties of the quantum doubles. The minimal condition that can be imposed for this goal is that both $W$ and $Z$ are second invertible as well
\cite{maj:fqgt}.

The results of the subsection \ref{wxxsolns},  can be reinterpreted as two--dimensional representations of the algebras
\be W_{12}T_1T_2=T_2T_1W_{12} \ee
for any invertible solution $W$ of the \YBE.

\section{Appendices}
\subsection{List of invertible solutions of the \YBE }\ll{apybeslns}
The important starting point for solving the WXZ system is the knowledge of complete set of invertible \sln s of the constant \YBE{} $[R,R,R]=0$. Up to the symmetries of the solution set of the \YBE
\be R'=\rho(T\ox T)R(T\ox T)^{-1},\ R'=R^+,\ R'=R^T,\ll{ybesyms}\ee
it  consists \cite{hie:ybecfn,hla:usybe,hla:baxtn} of the following (intersecting) subsets  parametrized by complex numbers $r,s,t,u,v$, 

Solution with three parameters
\[
R_{3.1}(r,s,t)  =  \left (\begin {array}{cccc} 1&0&0&0\\\noalign{\medskip}0&r&0&0
\\\noalign{\medskip}0&0&s&0\\\noalign{\medskip}0&0&0&t\end {array}
\right )  ,\  rst\neq 0 \]

Solutions with two parameters
\begin{eqnarray*}
 R_{2.1}(r,s) & = & \left (\begin {array}{cccc} 1&0&0&0\\\noalign{\medskip}0&r&0&0
\\\noalign{\medskip}0&1-rs&s&0\\\noalign{\medskip}0&0&0&1\end {array}
\right ) ,\  rs\neq 0
\\
 R_{2.2}(r,s) & = & \left (\begin {array}{cccc} 1&0&0&0\\\noalign{\medskip}0&r&0&0
\\\noalign{\medskip}0&1-rs&s&0\\\noalign{\medskip}0&0&0&-rs
\end {array}\right ) ,\  rs\neq 0
\\
R_{2.3}(u,v) & = & \left (\begin {array}{cccc} 1&0&0&0\\\noalign{\medskip}1&1&0
&0\\\noalign{\medskip}u&0&1&0\\\noalign{\medskip}v&u&1&1\end {array}
\right ) \end{eqnarray*} 

Solutions with one parameter
\begin{eqnarray*} 
 R_{1.1}(r) & = & \left (\begin {array}{cccc} r-{r}^{-1}+2&0&0&r-{r}^{-1}
\\\noalign{\medskip}0&r+{r}^{-1}&r-{r}^{-1}&0\\\noalign{\medskip}0&r-{
r}^{-1}&r+{r}^{-1}&0\\\noalign{\medskip}r-{r}^{-1}&0&0&r-{r}^{-1}-2
\end {array}\right ),\  r\neq 0
\\
 R_{1.2}(s) & = & \left (\begin {array}{cccc} 1&0&0&0\\\noalign{\medskip}0&
1&0&0\\\noalign{\medskip}0&1-s&s&0
\\\noalign{\medskip}1&0&0&-s\end {array}\right ),\  s\neq 0 \\
R_{1.3}(u) & = & \left (\begin {array}{cccc} 1&0&0&0\\\noalign{\medskip}-1&1&0
&0\\\noalign{\medskip}1&0&1&0\\\noalign{\medskip}u&-u&u&1
\end {array}\right ) 
\\
 R_{1.4}(t) & = & \left (\begin {array}{cccc} 0&0&0&1\\\noalign{\medskip}0&0&t&0
\\\noalign{\medskip}0&t&0&0\\\noalign{\medskip}1&0&0&0\end {array}
\right ),\  t\neq 0 \end{eqnarray*} 

Solutions without parameter
\begin{eqnarray*} R_{0.1} & = & \left (\begin {array}{cccc} 1&0&0&0\\\noalign{\medskip}0&1&0
&0\\\noalign{\medskip}0&0&1&0\\\noalign{\medskip}1&0&0&1\end {array}
\right ) \\
 R_{0.2} & = & \left (\begin {array}{cccc} 1&0&0&0\\\noalign{\medskip}0&-1&0
&0\\\noalign{\medskip}0&0&-1&0\\\noalign{\medskip}1&0&0&1\end {array}
\right ) \\
R_{0.3} & = & \left (\begin {array}{cccc} 1&0&0&i
\\\noalign{\medskip}0&1&1&0\\\noalign{\medskip}0&1&-1&0
\\\noalign{\medskip}i&0&0&1\end {array}\right ) \\
R_{0.4} & = & P = \left (\begin {array}{cccc} 1&0&0&0\\\noalign{\medskip}0&0&1
&0\\\noalign{\medskip}0&1&0&0\\\noalign{\medskip}0&0&0&1\end {array}
\right ) \\
\end{eqnarray*} 
\subsection{Equivalence classes of the $4\times4$ \mt ces with respect to the \tfn {} $ X'=(\unit\ox S)X(\unit\ox S)^{-1}$} \label{xansatz}
{\bf Lemma: } Any $4 \times 4$ matrix $A$ can be transformed into one of following 14 forms using transformation 
\begin{equation} A'= \lambda (\unit \ox S)A(\unit \ox  S)^{-1} 
\label{eqclassestransf}  \end{equation}
where $\lambda\in\complex,\ \unit = diag(1,1),\ S\in{\cal GL}(2,\complex)$.
\begin{eqnarray*}
A_{1} & = & \left ( \begin{array}{cccc}   a1 &   a2 &  \alpha+1 & 0
\\\noalign{\medskip}  b1 &   b2 & 0&  \alpha-1\\\noalign{\medskip}
 c1 &   c2 &   c3 &   c4\\\noalign{\medskip}  d1 &   d2
 &  d3 &   d4 \end{array} \right  )  \\
A_{2} & = & \left (\begin {array}{cccc}   a1 &   a2 &   a3 & 0
\\\noalign{\medskip}  b1 &   b2 & 1&  a3\\\noalign{\medskip} 
  c1 &   1 &   c3 &   c4\\\noalign{\medskip}0&  d2 &   
d3 &   d4 \end {array}\right  )   \\
A_{3} & = & \left (\begin {array}{cccc}   a1 &   a2 &   a3 & 0
\\\noalign{\medskip}  b1 &   b2 & 1&  a3\\\noalign{\medskip}
  \alpha+1 & 0&  c3 &   c4\\\noalign{\medskip}0&  \alpha-1 &   d3 & 
  d4 \end {array}\right  )   \\
A_{4} & = & \left (\begin {array}{cccc}   a1 &   a2 &   a3 & 0
\\\noalign{\medskip}  b1 &   b2 & 1&  a3\\\noalign{\medskip} 
  c1 & 0&  c3 &   c4\\\noalign{\medskip}  d1 &   c1 &   
d3 &   d4 \end {array}\right  )   \\
A_{5} & = & \left (\begin {array}{cccc}   a1 &   a2 &   a3 & 0
\\\noalign{\medskip}  b1 &   b2 & 0&  a3\\\noalign{\medskip} 
  \alpha+1 & 0&  c3 &   c4\\\noalign{\medskip}0&  \alpha-1 &   d3 & 
  d4 \end {array} \right  )  \\
A_{6} & = & \left (\begin {array}{cccc}   a1 &   a2 &   a3 & 0
\\\noalign{\medskip}  b1 &   b2 & 0&  a3\\\noalign{\medskip} 
  c1 & 0&  c3 &   1\\\noalign{\medskip}1&  c1 & 0&  d4
\end {array}\right  )    \\
A_{7} & = & \left (\begin {array}{cccc}   a1 &   a2 &   a3 & 0
\\\noalign{\medskip}  b1 &   b2 & 0&  a3\\\noalign{\medskip}
  c1 & 0&  \alpha+1 & 0\\\noalign{\medskip}1&  c1 & 0&  \alpha-1
\end {array}\right  )   \\
A_{8} & = & \left (\begin {array}{cccc}   a1 &   a2 &   a3 & 0
\\\noalign{\medskip}  b1 &   b2 & 0&  a3\\\noalign{\medskip} 
  c1 & 0&  c3 & 0\\\noalign{\medskip}1&  c1 &   d3 &   c3
\end {array}\right  )   \\
A_{9} & = & \left (\begin {array}{cccc}   a1 &   a2 &   a3 & 0
\\\noalign{\medskip}  b1 &   b2 & 0&  a3\\\noalign{\medskip} 
  c1 & 0&  \alpha+1 & 0\\\noalign{\medskip}0&  c1 & 0&  \alpha-1
\end {array}\right  )   \\
A_{10} & = & \left (\begin {array}{cccc}   a1 &   1 &   a3 & 0
\\\noalign{\medskip}0&  b2 & 0&  a3\\\noalign{\medskip}  c1 & 0
&  c3 & 0\\\noalign{\medskip}0&  c1 & 1&  c3 \end {array}
\right  )   \\
A_{11} & = & \left (\begin {array}{cccc}   \alpha+1 & 0&  a3 & 0
\\\noalign{\medskip}0&  \alpha-1
 & 0&  a3\\\noalign{\medskip}  c1 & 0
&  c3 & 0\\\noalign{\medskip}0&  c1 & 1&  c3 \end {array}
\right  )  \\ 
A_{12} & = & \left (\begin {array}{cccc}   a1 & 0&  a3 & 0
\\\noalign{\medskip}  b1 &   a1 & 0&  a3\\\noalign{\medskip} 
  c1 & 0&  c3 & 0\\\noalign{\medskip}0&  c1 & 1&  c3
\end {array}\right  )  \\
A_{13} & = & \left (\begin {array}{cccc}   a1 & 0&  a3 & 0
\\\noalign{\medskip}0&  b2 & 0&  a3\\\noalign{\medskip}  c1 & 0
&  c3 & 0\\\noalign{\medskip}0&  c1 & 0&  c3 \end {array}
\right  )   \\
A_{14} & = & \left (\begin {array}{cccc}   a1 & 0&  a3 & 0
\\\noalign{\medskip}1&  a1 & 0&  a3\\\noalign{\medskip}  c1 & 0
&  c3 & 0\\\noalign{\medskip}0&  c1 & 0&  c3 \end {array}
\right  )   \\
\end{eqnarray*}
{\bf Proof: } In order to prove previous lemma it is convenient to write
matrix A divided into square blocks, i.e. 
\begin{equation} A = \left ( \begin{array}{cc} A_{11} & A_{12} \\ A_{21} & A_{22} \end{array} \right ). \label{eqclassestransf1}
\end{equation}
Blocks $A_{11},\ldots,A_{22}$ are transforming like $A_{ij}'= \lambda S A_{ij} S^{-1}$ under tranformation (\ref{eqclassestransf}). It is possible to convert one chosen block into Jordan's canonical form, i.e. one of the following forms:
\begin{equation} 
1) \left ( \begin{array}{cc} \alpha & 0 \\ 0 & \alpha \end{array} \right ) \qquad 
2) \left ( \begin{array}{cc} \alpha & 0 \\ 1 & \alpha \end{array} \right ) \qquad
3) \left ( \begin{array}{cc} \alpha+1 & 0 \\ 0 & \alpha-1 \end{array} \right )
\end{equation}
(~see Jordan's theorem in linear algebra~).
Case $1)$ is of course invariant w.r.t. (\ref{eqclassestransf}), in case $2)$ it is possible to convert  another arbitrary block into upper triangular (~with non-diagonal element equal 1~) or lower triangular (~with identical diagonal elements~) or diagonal block. Using these properties, one can find given list of classes, firstly simplifing upper right block, secondly using remaining symmetries to simplify lower left block, then lower right block and finally upper left block.

\subsection{List of  $X$ \mt ces}\label{xseznam}
Solutions of the equation \rf{wxx1}) can be written in terms of the following matrices.
\newcounter{xmatice}
\begin{eqnarray*}
\refstepcounter{xmatice}\label{diag}
X_{\thexmatice}(a,b,c) & = & \left ( \begin{array}{cccc} 
  1 & 0 & 0 & 0 \\
  0 & a & 0 & 0 \\
  0 & 0 & b & 0 \\
  0 & 0 & 0 & c
  \end{array} \right )
\\
\refstepcounter{xmatice}\label{X.jj1}
X_{\thexmatice}(a,b,c) & = & \left ( \begin{array}{cccc} 
  1 & 0 & 0 & 0 \\
  b & 1 & 0 & 0 \\
  0 & 0 & a & 0 \\
  0 & 0 & c & a
  \end{array} \right ),\ \  b \neq 0\  {\rm or}\ c\neq 0
\\
\refstepcounter{xmatice}\label{X.jj2}
X_{\thexmatice}(a,b,c) & = & \left ( \begin{array}{cccc} 
  1 & 0 & 0 & 0 \\
  0 & a & 0 & 0 \\
  b & 0 & 1 & 0 \\
  0 & c & 0 & a
  \end{array} \right ),\ 
\\
\refstepcounter{xmatice}\label{X10}
X_{\thexmatice}(a,b,c) & = & \left ( \begin{array}{cccc} 
  1 & 0 & 0 & 0 \\
  a & 1 & 0 & 0 \\
  b & 0 & 1 & 0 \\
  c & b & a & 1
  \end{array} \right ),\ \  a \neq 0\  {\rm or}\ b \neq 0\  {\rm or}\ c\neq 0
\\
\refstepcounter{xmatice}\label{X8.3}
X_{\thexmatice}(a,b,c) & = & \left ( \begin{array} {cccc}
  1 & 0 & 0 & 0 \\
  0 & 0 & 0 & a \\
  0 & 0 & b & 0 \\
  0 & c & 0 & 0
  \end{array} \right ),\ \  abc \neq 0 
\\
\refstepcounter{xmatice}\label{X8.4}
X_{\thexmatice}(a,b,c) & = & \left ( \begin{array} {cccc}
  0 & 0 & 1 & 0 \\
  0 & 0 & 0 & a \\
  b & 0 & 0 & 0 \\
  0 & c & 0 & 0
  \end{array} \right ), \ \  abc \neq 0 
\\
\refstepcounter{xmatice}\label{X8.5}
X_{\thexmatice}(a,b) & = & \left ( \begin{array} {cccc}
  0 & 0 & 1 & 0 \\
  0 & 0 & a & 1 \\
  1 & 0 & 0 & 0 \\
  b & 1 & 0 & 0
  \end{array} \right ) ,\ \  a \neq 0\  {\rm or}\ b\neq 0
\\
\refstepcounter{xmatice}\label{X8.7}
X_{\thexmatice}(a) & = & \left ( \begin{array} {cccc}
  0 & 0 & 1 & 0 \\
  0 & 0 & 0 & -1 \\
  0 & a & 0 & 0 \\
  1 & 0 & 0 & 0
  \end{array} \right ) 
\\
\refstepcounter{xmatice}\label{h5}
X_{\thexmatice}(b,c,d) & = & \left ( \begin{array} {cccc}
  1 & 0 & 0 & 0 \\
  0 & b & 0 & 0 \\
  0 & 0 & c & 0 \\
  1 & 0 & 0 & d
  \end{array} \right )
\\   
\refstepcounter{xmatice}\label{h6}
X_{\thexmatice}(b,c,d,g) & = & \left ( \begin{array} {cccc}
  1 & 0 & 0 & 0 \\
  0 & b & 1 & 0 \\
  0 & 0 & c & 0 \\
  g & 0 & 0 & d
  \end{array} \right ) \ g\neq 0
\\
\refstepcounter{xmatice}\label{X.30}
X_{\thexmatice}(a,b,c,d) & = & \left ( \begin{array}{cccc} 
  0 & ab & 1 & 0 \\
  a & 0 & 0 & -1 \\
  d & 0 & 0 & b \\
  0 &-d & 1 & 0
  \end{array} \right ),
\\  
\refstepcounter{xmatice}\label{h7}
X_{\thexmatice}(a,b) & = & \left ( \begin{array}{cccc} 
  a & 0 & 0 & 0 \\
  0 & a & 0 & 0 \\
  1 & 0 & 0 & b \\
  0 &-1 & b & 0
  \end{array} \right )
\\  
\refstepcounter{xmatice}\label{octag}
X_{\thexmatice}(g) & = & \left ( \begin{array}{cccc} 
  0 & 0 & i+1 & 0 \\
  2(i+1)g & 0 & 0 & i-1 \\
  -ig & 0 & 0 & 1 \\
  0 & g & 0 & 0
  \end{array} \right )
\\
\refstepcounter{xmatice}\label{X2.1}
X_{\thexmatice}(a,\epsilon,p) &=& \left ( \begin{array}{cccc}
\epsilon ap & \epsilon c & 1 & 0 \\
     \epsilon c& \epsilon a p & 0 & -1 \\
         \epsilon & 0 & a & cp  \\
         0 & -\epsilon & cp & a
       \end{array} \right ) ,\ c^2p={a^2p-1}
\\
\refstepcounter{xmatice}\label{X9.9}
X_{\thexmatice}(a) & = & \left ( \begin{array} {cccc}
  a & 1-a^2 & 0 & 0 \\
  1 & -a & 0 & 0 \\
  0 & 0 & 1 & 0 \\
  0 & 0 & 0 & -1
  \end{array} \right ) ,\ \ a\neq\pm 1
\\
\refstepcounter{xmatice}\label{h11.2}
X_{\thexmatice}(a,b,c,d) &=& \left ( \begin{array}{cccc}
         1 & 0 & 0 & 0 \\
         a & 1 & 0 & 0 \\
         b & 0 & 1 & 0 \\
         0 & c & d & 1
       \end{array} \right ),\ a\neq d,b\neq c 
\\
\refstepcounter{xmatice}\label{h11.1}
X_{\thexmatice}(a,b) &=& \left ( \begin{array}{cccc}
         1 & 0 & 0 & 0 \\
         a & 1 & 1 & 0 \\
         a+b+1 & 0 & 1 & 0 \\
         0 & a+b-1 & b & 1
       \end{array} \right ) 
\\
\refstepcounter{xmatice}\label{X9.4}
X_{\thexmatice}(a) & = & \left ( \begin{array} {cccc}
  0 & 0 & 1 & 0 \\
  0 & 0 & 0 & -1 \\
  a & 1 & 0 & 0 \\
  1-a^2 & -a & 0 & 0
  \end{array} \right ) 
\\
\refstepcounter{xmatice}\label{X10.5}
X_{\thexmatice}(a,b,c,\epsilon) & = & \left ( \begin{array}{cccc}
  1 & 0 & 0 & 0 \\
  0 & a & 0 & 0 \\
  b & 0 & \epsilon & 0 \\
  0 & c & 0 & -a
  \end{array} \right )
\\
\refstepcounter{xmatice}\label{X10.6}
X_{\thexmatice}(a,b,c) & = & \left ( \begin{array}{cccc} 
  1 & 0 & 0 & 0 \\
  b & 1 & 0 & 0 \\
  c & 0 & -1 & 0 \\
  a & c & -b & -1
  \end{array} \right )
\\
\refstepcounter{xmatice}\label{X1.4}
X_{\thexmatice}(a,b) & = & \left ( \begin{array} {cccc}
  a & b & 1 & 0 \\
  b & a & 0 & -1 \\
  0 & i & -b & -a \\
  -i & 0 & -a & -b
  \end{array} \right ) 
\\
\refstepcounter{xmatice}\label{X8.15}
X_{\thexmatice}(a,b,c,d,e,f,g,h) & = &  \left ( \begin{array}{cccc}
         a & 0 & b & 0 \\
         c & a & d & b \\
         e & 0 & f & 0 \\
         g & e & h & f
       \end{array} \right ) 
\\
\refstepcounter{xmatice}\label{X8.16}
X_{\thexmatice}(a,b,c,d,e,f,g,h) &=& \left ( \begin{array}{cccc}
         a & 0 & b & 0 \\
         0 & c & 0 & d \\
         e & 0 & f & 0 \\
         0 & g & 0 & h
       \end{array} \right ) 
\end{eqnarray*}

\end{document}